\documentclass[12pt,a4paper,reqno]{amsart}
\usepackage{amsmath}
\usepackage{amsfonts}
\usepackage{amssymb}
\numberwithin{equation}{section}

\theoremstyle{plain}

\newtheorem{theorem}[subsection]{Theorem}
\newtheorem{proposition}[subsection]{Proposition}
\newtheorem{lemma}[subsection]{Lemma}
\newtheorem{corollary}[subsection]{Corollary}
\newtheorem{conjecture}[subsection]{Conjecture}
\newtheorem{observation}[subsection]{Observation}

\theoremstyle{definition}

\newtheorem{definition}[subsection]{Definition}

\newtheorem{question}[subsection]{Question}

\renewcommand{\leq}{\leqslant}
\renewcommand{\geq}{\geqslant}

\newsavebox{\proofbox}
\savebox{\proofbox}{\begin{picture}(7,7)%
  \put(0,0){\framebox(7,7){}}\end{picture}}
\def\boxeq{\tag*{\usebox{\proofbox}}}

\newcommand\codim{\operatorname{codim}}
\newcommand\rk{\operatorname{rk}}
\newcommand\Spec{\operatorname{Spec}}
\def\proof{\noindent\textit{Proof. }}
\def\endproof{\hfill{\usebox{\proofbox}}}

\def\E{\mathbb{E}}
\def\Z{\mathbb{Z}}

\def\R{\mathbb{R}}
\def\F{\mathbb{F}}
\def\f{\mathbf{f}}
\def\B{\mathcal{B}}
\def\C{\mathbb{C}}
\def\b{\mathbf{b}}
\def\GL{\operatorname{\mathfrak{M}}}

\def\sym{\operatorname{sym}}

\begin{document}
\title{Montr\'eal notes on quadratic fourier analysis}
\author{Ben Green}
\address{Centre for Mathematical Sciences\\ Wilberforce Road\\ Cambridge CB3 0WA\\ England}
\email{b.j.green@dpmms.cam.ac.uk}

\thanks{The author is a Clay Research fellow. He is grateful for the support of the Clay Mathematics Institute, which enabled him to attend the activities in Montr\'eal.}

\begin{abstract}
These are notes to accompany four lectures that I gave at the \emph{School on additive combinatorics,} held in Montr\'eal, Qu\'ebec between March 30th and April 5th 2006. 

My aim is to introduce ``quadratic fourier analysis'' in so far as we understand it at the present time. Specifically, we will describe ``quadratic objects'' of various types and their relation to additive structures, particularly four-term arithmetic progressions.

I will focus on qualitative results, referring the reader to the literature for the many interesting quantitative questions in this theory. Thus these lectures have a distinctly ``soft'' flavour in many places. 

Some of the notes cover unpublished work which is joint with Terence Tao. This will be published more formally at some future juncture.
\end{abstract}

\maketitle

\section{Lecture 1}

Topics to be covered:

\begin{itemize}
\item Introduction. The finite field philosophy.
\item Review of notation and basic properties of the Fourier transform
\item Counting 3- and 4-term arithmetic progressions using the Gowers $U^2$- and $U^3$-norms: generalised von Neumann theorems.
\item Inverse theorem for the Gowers $U^2$-norm.
\item The ``quadratic'' example for the Gowers $U^3$-norm. 
\item Brief revision of key results from additive combinatorics.
\end{itemize}

What is ``quadratic Fourier analysis?''. The aim of this series of lectures is to give a reasonably detailed answer to that question, at least in so far as is possible at the present time. 

It would, however, be presumptuous to suppose that any reader would venture to the end of these notes in order to discover the meaning of the title, so we begin with a very brief introduction.

Fourier analysis, or ``linear'' Fourier analysis as we shall call it in these notes, is a multi-faceted subject. One rather small part of it is concerned with solving linear equations. Two examples of theorems which may be proven using some kind of study of the Fourier transform are
\begin{itemize}
\item (Chowla/van der Corput) There are infinitely many 3-term arithmetic progressions of primes.
\item (Roth) Let $\delta > 0$ be fixed. Then if $N > N_0(\delta)$ is sufficiently large, any subset $A \subseteq \{1,\dots,N\}$ with size at least $\delta N$ contains three distinct elements in arithmetic progression.
\end{itemize}
Note that an arithmetic progression of length three is defined by a single linear equation $x_1 + x_3 = 2x_2$.

Standard Fourier analysis fails in many situations where we are interested in a pair of linear equations. The natural example here is a progression of length four, which is defined by the equations $x_1 + x_3 = 2x_2$, $x_2 + x_4 = 2x_3$. This is the situation where quadratic Fourier analysis is appropriate. Thus by developing the methods that we will talk about in these lectures, it is possible to prove
\begin{itemize}
\item (Green-Tao) There are infinitely many 4-term arithmetic progressions of primes.
\item (Szemer\'edi) Let $\delta > 0$ be fixed. Then if $N > N_0(\delta)$ is sufficiently large, any subset $A \subseteq \{1,\dots,N\}$ with size at least $\delta N$ contains three distinct elements in arithmetic progression.
\end{itemize}
In fact we will not prove either of these theorems in this course, since we will be working in a model setting. A common theme in additive combinatorics is the consideration of \emph{finite field models}. A full discussion may be found in  \cite{green-fin-field}, but the basic idea is as follows. For many problems in additive combinatorics one is interested in the interval $\{1,\dots,N\}$. However, it is convenient to work in a group, and so one often uses various technical devices in order to place the problem at hand in $\Z/N\Z$. Once this is done, it is often easy to formulate an analogous question inside an arbitrary finite abelian group $G$. In most applications that we know of, this more general problem is scarcely harder to solve than in the specific case $G = \Z/N\Z$. However, there is a family of groups, namely the groups $\F_p^n$ where $p$ is a small prime, in which it can be relatively easy to work. Techniques used to prove theorems in this setting can often be used to guide proof techniques in $\Z/N\Z$, which provide theorems of actual number theoretic interest.

In this series of lectures we will focus almost exclusively on the group $G = \F_5^n$. I am rather fond of the prime $5$ since it is the smallest for which the notion of a 4-term arithmetic progression is sensible.

We will conclude with a discussion of the general case at the end, in as much detail as time permits. It turns out that the theory for $\Z/N\Z$ is surprisingly rich, and there are strong connections with the ergodic theory techniques that are discussed in the lectures of Bryna Kra in these Proceedings. 

\emph{Notation.} Opinion seems to be converging in additive combinatorics about what constitutes the ``standard'' notation, and I will endeavour to keep to these norms. If $X$ is any finite set and $f : X \rightarrow \C$ is any function then we write
\[ \E_{x \in X} f(x) := |X|^{-1} \sum_{x \in X} f(x).\]
This means that it is often possible to avoid worrying about normalising factors. 

Unless specified otherwise, we will set $G := \F_5^n$ and write $N := |G| = 5^n$.
Any character on $G$ (that is, homomorphism $\gamma : G \rightarrow \C^{\times}$) has the form $x \mapsto \omega^{r^T x}$, where $\omega := e^{2\pi i /5}$ and $r \in \F_5^n$ is a vector. We write $\widehat{G}$ for $\F_5^n$ when considered as the group of characters in this way.

If $f : G \rightarrow \C$ is any function then we define its Fourier transform $\widehat{f} : \widehat{G} \rightarrow \C$ by
\[ \widehat{f}(r) := \E_{x \in G} f(x) \omega^{r^T x}.\] We distinguish the \emph{trivial character} corresponding to $r = 0$, which takes the value $1$ for all $x \in G$. If $f,g : G \rightarrow \C$ are two functions then we define the \emph{convolution} $f \ast g : G \rightarrow \C$ by
\[ (f \ast g)(x) := \E_{y \in G} f(x)g(y - x).\]
Note that when working on $G$ we always use the \emph{Haar measure} which assigns weight $|G|^{-1}$ to any $x \in G$. When working on $\widehat{G}$ we use the \emph{counting measure} which assigns weight $1$ to every $r \in \widehat{G}$. These measures are dual to one another, which in practice means that in formul{\ae} such as those in Lemma \ref{basic-props} below one can simply write $\E_{x \in G}$ and $\sum_{r \in \widehat{G}}$, and thereafter be untroubled by normalising factors.

When we talk of $L^p$ norms, these will always be taken with respect to the appropriate underlying measure. Thus
\[ \Vert f \Vert_1 := \E_{x \in G} |f(x)|,\] whereas
\[ \Vert \widehat{f} \Vert_4 := \big( \sum_{r \in \widehat{G}} |\widehat{f}(r)|^4 \big)^{1/4}.\]

I will be assuming familiarity with the basic properties of the Fourier transform, which are all straightforward consequences of the orthogonality relations
\[ \sum_{r} \omega^{r^T x} = \left\{\begin{array}{ll} N & \mbox{if $x = 0$} \\ 0 & \mbox{otherwise} \end{array} \right.\]
and 
\[ \E_{x \in G} \omega^{r^T x} = \left\{\begin{array}{ll} 1 & \mbox{if $r = 0$} \\ 0 & \mbox{otherwise.} \end{array} \right.\]
\begin{lemma}[Basic properties of the Fourier transform]\label{basic-props}
Suppose that $f,g : G \rightarrow \C$ are any two functions. Then
\begin{enumerate}
\item We have $\widehat{f}(0) = \E_{x \in G} f(x)$. For any $r$ we have $|\widehat{f}(r)| \leq \Vert f \Vert_1$.
\item \textup{(Parseval identity)} We have \[ \E_{x \in G} f(x) \overline{g}(x) = \sum_{r \in \widehat{G}} \widehat{f}(r) \overline{\widehat{g}(r)}.\] In particular $\Vert f \Vert_2 = \Vert \widehat{f} \Vert_2$.
\item \textup{(Inversion)} We have
\[ f(x) = \sum_{r \in \widehat{G}} \widehat{f}(r) \omega^{-r^T x}.\]
\item \textup{(Convolution)} We have $(f \ast g)^{\wedge} = \widehat{f} \widehat{g}$.
\end{enumerate}
\end{lemma}
The last item here illustrates how we will denote the Fourier transform of expressions $E$ for which it would be too cumbersome to write $\widehat{E}$.

Let us now start with the main business of these lectures.
Let $G$ be a finite abelian group with order $N$ which is coprime to $6$, and let $f_1,\dots,f_4 : G \rightarrow [-1,1]$ be functions. In these notes a central r\^ole will be played by the multilinear operators $\Lambda_3$ and $\Lambda_4$, defined by
\[ \Lambda_3(f_1,f_2,f_3) := \E_{x,d} f_1(x)f_2(x+d)f_3(x+2d)\]
and
\[ \Lambda_4(f_1,f_2,f_3,f_4) := \E_{x,d} f_1(x)f_2(x+d)f_3(x+2d)f_4(x+3d).\]
Thus $\Lambda_3$ counts the number of 3-term arithmetic progressions ``along the $f_i$'', whilst $\Lambda_4$ counts the number of 4-term progressions.\footnote{Whilst we will talk exclusively about 3- and 4-term arithmetic progressions, the reader should note that much of what we have to say may be adapted to more general problems where it is of interest to count the number of solutions to a linear equation, or to a pair of linear equations.} 

When the functions $f_i$ are characteristic functions, the operators $\Lambda_3$ and $\Lambda_4$ may be interpreted combinatorially.

\begin{observation}
Suppose that $f_i = 1_{A_i}$, where $A_i \subseteq G$ is a set. Then $\Lambda_3(1_{A_1},1_{A_2},1_{A_3})$ is equal to $N^{-2}$ times the number of triples $(a_1,a_2,a_3) \in A_1 \times A_2 \times A_3$ which are in arithmetic progression. Similarly, $\Lambda_4(1_{A_1},1_{A_2},1_{A_3},1_{A_4})$ is equal to $N^{-2}$ times the number of quadruples $(a_1,a_2,a_3,a_4) \in A_1 \times A_2 \times A_3 \times A_4$ which are in arithmetic progression.
\end{observation}

There are certainly many situations in which one might be interested in counting the number of 3- or 4-term progressions inside a set. To do this, we normally proceed as follows. If $A \subseteq G$ is a set with size $\alpha N$, then write $f_A := 1_A - \alpha$. This is called the \emph{balanced function} of $A$, and it has expected value 0. 

\begin{lemma}[Balanced function decomposition]\label{bal-decomp}
Suppose that $A_1,\dots,A_4 \subseteq G$, and that $|A_i| = \alpha_i N$. Then we have 
\begin{equation}\label{seven} \Lambda_3(1_{A_1},1_{A_2},1_{A_3}) = \alpha_1\alpha_2\alpha_3 + \mbox{\textup{(seven other terms)}},\end{equation}
where each of the seven terms has the form $\Lambda_3(g_1,g_2,g_3)$ where each $g_i$ is either $f_{A_i}$ or $\alpha_i$ and at least one is equal to $f_{A_i}$. Similarly
\begin{equation}\label{fifteen} \Lambda_4(1_{A_1},1_{A_2},1_{A_3},1_{A_4}) = \alpha_1\alpha_2\alpha_3\alpha_4 + \mbox{\textup{(fifteen other terms)}},\end{equation}
where each of the fifteen terms has the form $\Lambda_4(g_1,g_2,g_3,g_4)$ where each $g_i$ is either $f_{A_i}$ or $\alpha_i$ and at least one is equal to $f_{A_i}$.
\end{lemma}

Let us specialise to the case $A_1 = A_2 = A_3 = A_4 = A$ for simplicity, and write $|A| = \alpha N$. What do we ``expect'' $\Lambda_3(1_A,1_A,1_A)$ and $\Lambda_4(1_A,1_A,1_A,1_A)$ to be? It is not hard to see that for a ``random'' set $A$, generated by tossing a coin which comes up heads with probability $\alpha$ to decide whether each $x \in G$ lies in $A$, the expected value of $\Lambda_3(1_A,1_A,1_A)$ is approximately $\alpha^3$, whilst the expected value of $\Lambda_4(1_A,1_A,1_A,1_A)$ is approximately $\alpha^4$. Note that these quantities are exactly the ``main terms'' in the expansions of Lemma \ref{bal-decomp}. It is thus reasonable to suggest that the other seven terms in \eqref{seven} measure some kind of ``non-uniformity'' of $A$ relevant to 3-term progressions, whilst the fifteen terms in \eqref{fifteen} do the same for 4-term progressions.

Let us make a preliminary definition.

\begin{definition}[Uniformity along progressions]\label{unif-ap-def}
Let $A \subseteq G$ be a set with $|A| = \alpha N$, and let $f_A := 1_A - \alpha$ be the balanced function of $A$. Let $\delta \in (0,1)$ be a parameter. Then we say that $A$ exhibits $\delta$-uniformity along 3-term progressions if whenever we have three functions $g_1,g_2,g_3 \rightarrow [-1,1]$, at least one of which is equal to $f_A$, then
\[ |\Lambda_3(g_1,g_2,g_3)| \leq \delta.\]
We define non-unifomity along 4-term progressions similarly.
\end{definition}

\noindent\textit{Remark.} It is not, at first sight, obvious that there are \emph{any} sets which are uniform along progressions. 

\begin{lemma}
Suppose that $A \subseteq G$ is a set with $|A| = \alpha N$. If $A$ is $\delta$-uniform along 3-term progressions, then
\[ |\Lambda_3(1_A,1_A,1_A) - \alpha^3| \leq 7\delta.\]
If $A$ is $\delta$-uniform along $4$-term progressions, then
\[ |\Lambda_4(1_A,1_A,1_A,1_A) - \alpha^4| \leq 15\delta.\]
\end{lemma}
\proof Immediate consequence of Lemma \ref{bal-decomp}.\endproof

The following question will be a recurring theme of these lectures:

\begin{question}\label{main-question}
Suppose that $A$ is not $\delta$-uniform along 3- or 4-term progressions. Can we say something ``useful'' about $A$?
\end{question}

Of course, the notion of ``useful'' is a subjective one. The reader may assume, however, that the mere failure of Definition \ref{unif-ap-def} does not constitute ``useful''. We will see that if $A$ is not uniform along 3-term progressions, then it exhibits ``linear'' behaviour, whilst functions which are not uniform along 4-term progressions are somehow ``quadratic''. 

Formulating, proving, and using statements of this type is our main goal in these notes.

Question \ref{main-question} may be answered very satisfactorily using Fourier analysis. The key tool is the following simple lemma, whose proof is an amusing exercise using the basic properties of the Fourier transform.

\begin{lemma}\label{four-form}
Let $f_1,f_2,f_3 : G \rightarrow \R$ be any three functions. Then
\begin{equation}\boxeq \Lambda_3(f_1,f_2,f_3) = \sum_r \widehat{f}_1(r)\widehat{f}_2(-2r)\widehat{f}_3(r)\end{equation}
\end{lemma}

\begin{proposition}[Inverse result for 3-term progressions, I]\label{3-inv-1}
Suppose that $A$ is not $\delta$-uniform along 3-term progressions. Then $\Vert \widehat{f}_A \Vert_{\infty} \geq \delta$, that is to say there is some $r \in \widehat{G}$ such that $|\widehat{f}_A(r)| \geq \delta$.
\end{proposition}
\proof Suppose that 
\[ |\Lambda_3(g_1,g_2,f_A)| \geq \delta\] for some functions $g_1,g_2 : G \rightarrow [-1,1]$ (the analysis of the other two cases, when $g_1 = f_A$ or $g_2 = f_A$, is more-or-less identical).
We have, by Lemma \ref{four-form} the formula
\[ \Lambda_3(g_1,g_2,f_A) = \sum_{r} \widehat{g}_1(r)\widehat{g}_2(-2r)\widehat{f_A}(r).\]
Thus by Cauchy-Schwarz and Parseval's identity we infer that
\begin{equation}\boxeq \delta \leq |\sum_{r} \widehat{g}_1(r)\widehat{g}_2(-2r)\widehat{f_A}(r)| \leq \Vert \widehat{f}_A \Vert_{\infty} \Vert \widehat{g_1} \Vert_2 \Vert \widehat{g_2} \Vert_2 \leq \Vert \widehat{f}_A \Vert_{\infty}.\end{equation}

This is a very clean result, but the method of proof (appealing to a formula in Fourier analysis) has not, so far, proved amenable to generalisation. One way to generalise an argument is to first try and find a more longwinded, less natural looking approach and try and generalise that. We will describe such an approach now, though we hope that any reader looking back on this section later on will not consider it so unnatural. Note that the result is the same as Proposition \ref{3-inv-1}, but the bound is slightly worse.

\begin{proposition}[Inverse result for 3-term progressions, II]\label{3-inv-2}
Suppose that $A$ is not $\delta$-uniform along 3-term arithmetic progressions. Then $\Vert \widehat{f}_A \Vert_{\infty} \geq \delta^2$.
\end{proposition}
\proof Let us first observe that 
\[ \Lambda_3(g_1,g_2,f_A) = \E_{y_1,y_2} g_1(-y_1)g_2(\textstyle\frac{1}{2} y_2\displaystyle)f_A(y_1 + y_2).\] This is a simple reparametrisation. Applying the Cauchy-Schwarz inequality, we have
\begin{align*} |\Lambda_3(g_1,g_2,f_A)|^2 &\leq \E_{y_2} |\E_{y_1} g_1(-y_1) f_A(y_1 + y_2)|^2 \\ & = \E_{y_1,y'_1,y_2} f_A(y_1 + y_2)f_A(y'_1 + y_2)g_1(-y_1)g_1(-y'_1).\end{align*}
Applying Cauchy-Schwarz again, we have
\begin{align}\nonumber |\Lambda_3(g_1,g_2,f_A)|^4 &\leq \E_{y_1,y'_1}|\E_{y_2} f_A(y_1 + y_2)f_A(y'_1 + y_2)|^2 \\ &= \E_{y_1,y'_1,y_2,y'_2} f_A(y_1 + y_2)f_A(y'_1 + y_2)f_A(y_1 + y'_2)f_A(y'_1 + y'_2). \label{in-a-min} \end{align}

This last expression is called the (fourth power of) the \emph{Gowers $U^2$-norm of $f_A$}. Thus we define
\begin{equation}\label{gowers-u2-def} \Vert f_A \Vert_{U^2}^4 := \E_{y_1,y'_1,y_2,y'_2} f_A(y_1 + y_2)f_A(y'_1 + y_2)f_A(y_1 + y'_2)f_A(y'_1 + y'_2).\end{equation}
It is often useful to write this in the alternative form
\[ \Vert f_A \Vert_{U^2}^4 = \E_{x,h_1,h_2} f_A(x)f_A(x+h_1)f_A(x+h_2)f_A(x+h_1+h_2).\]
It is not hard to show that $\Vert \cdot \Vert_{U^2}$ is a norm using the Cauchy-Schwarz inequality several times. We will not make much use of this fact, and refer the reader to \cite{gowers-longaps} for the proof.

Note that \eqref{in-a-min} implies that if $|\Lambda_3(g_1,g_2,f_A)| \geq \delta$ then $\Vert f_A \Vert_{U^2} \geq \delta$.

What now? Another way to see that $\Vert \cdot \Vert_{U^2}$ is a norm is to observe that
\[ \Vert f \Vert_{U^2}^4 = \Vert f \ast f \Vert_2^2 = \Vert (f \ast f)^{\wedge} \Vert_2^2 = \Vert \widehat{f} \Vert_4.\]
Thus if $\Vert f_A \Vert_{U^2} \geq \delta$ then we have
\[ \delta^4 \leq \Vert \widehat{f}_A \Vert_4 \leq \Vert \widehat{f}_A \Vert_{\infty}^2 \Vert \widehat{f}_A \Vert_2^2 \leq \Vert \widehat{f}_A \Vert_{\infty}^2,\]
which concludes the proof in the case that $|\Lambda_3(g_1,g_2,f_A)| \geq \delta$. Again, the cases when $f_A = g_1$ or $g_2$ can be dealt with very similarly, and are left to the reader; the parametrisations leading to \eqref{in-a-min} must be modified slightly.\endproof

At the moment, it is hard to see what has been gained here. To prove the result, we still had to fall back on a formula of Fourier analysis, and furthermore the bound we obtain is worse than that in Proposition \ref{3-inv-1}.

We may summarise the argument in Proposition \ref{3-inv-2} as follows, giving the two distinct parts a name.

\begin{itemize}
\item (Generalised von Neumann theorem) The operator $\Lambda_3$ is controlled by the Gowers $U^2$-norm. Specifically for any three functions $f_1,f_2,f_3 : G \rightarrow [-1,1]$ we have
\[ |\Lambda_3(f_1,f_2,f_3)| \leq \inf_{i=1,2,3} \Vert f_i \Vert_{U^2}.\]
\item (Gowers inverse theorem) If the Gowers $U^2$-norm of a function $f : G \rightarrow [-1,1]$ is large, $f$ must have a large Fourier coefficient:
\[ \Vert f \Vert_{U^2} \geq \delta \quad \Rightarrow \quad \Vert \widehat{f} \Vert_{\infty} \geq \delta^2.\]
\end{itemize}

We note that the Gowers inverse theorem is necessary and sufficient. Indeed if $\Vert \widehat{f} \Vert_{\infty} \geq \delta$ then clearly $\Vert \widehat{f} \Vert_4 \geq \delta$, and so of course $\Vert f \Vert_{U^2} \geq \delta$.

This division of labour into two parts turns out to be the natural way to proceed for $\Lambda_4$ (and higher operators). The first part of the argument (the definition of the Gowers norm and the Generalised von Neumann theorem) goes through somewhat straightforwardly. The second part (the Gowers inverse theorem) does not, since we do not know of a formula analagous to $\Vert f \Vert_{U^2} = \Vert \widehat{f} \Vert_4$. 

\begin{definition}[Gowers $U^3$-norm]
Let $f : G \rightarrow [-1,1]$ be a function. Then we define
\begin{align*} \Vert f \Vert_{U^3}^8 := & \E_{\substack{y_1,y_2,y_3\\y'_1,y'_2,y'_3}} f(y_1 + y_2 + y_3)f(y'_1 + y_2 + y_3)f(y_1 + y'_2 + y_3)f(y_1 + y_2 + y'_3) \times \\ & \qquad\qquad\times f(y'_1 + y'_2 + y_3)f(y'_1 + y_2 + y'_3)f(y_1 + y'_2 + y'_3)f(y'_1 + y'_2 + y'_3) \\& = \E_{x,h_1,h_2,h_3 \in G} f(x)f(x+h_1) f(x+h_2) \times \\ & \times f(x+h_3)f(x+h_1+h_2)f(x+h_1 + h_3)f(x+h_2 + h_3) f(x+h_1+h_2 + h_3).\end{align*}
\end{definition}
Note that this is a kind of sum of $f$ over 3-dimensional parallelepipeds. We omit the proof that $\Vert f \Vert_{U^3}$ is actually a norm (see \cite{gowers-longaps}).

\begin{proposition}[Generalised von Neumann theorem for 4-term APs]\label{4ap-gvn} Let $f_1,\dots,f_4 : G \rightarrow [-1,1]$ be any four functions. Then we have
\[ |\Lambda_4(f_1,\dots,f_4)| \leq \inf_{i = 1,\dots,4} \Vert f_i \Vert_{U^3}.\]
In particular if $A$ is not $\delta$-uniform along four-term progressions then $\Vert f_A \Vert_{U^3} \geq \delta$.
\end{proposition}
\proof The idea is the same as in Proposition \ref{3-inv-2}. Here, we find a suitable reparametisation of $\Lambda_4(f_1,\dots,f_4)$, and then apply the Cauchy-Schwarz inequality three times. A ``suitable reparametrisation'' turns out to be
\begin{align}\nonumber\Lambda_4 &(f_1,f_2,f_3,f_4)\\ & = \E_{y_1, y_2, y_3 \in G} \textstyle f_1(-\frac{1}{2} y_2 - 2y_3) f_2(\frac{1}{3} y_1 - y_3) f_3(\frac{2}{3} y_1 + \frac{1}{2}y_2) f_4(y_1 + y_2 + y_3).\label{nice-param}\end{align}
For the rest of this section let $\b()$ denote any function bounded by $1$. Different occurrences of $\b$ may denote different functions. The Cauchy-Schwarz inequality implies that
\begin{equation}\label{cz-repeat} |\E_{x \in X}\E_{y \in Y} \b(x) f(x,y)| \leq |\E_{x \in X} \E_{y^{(0)}, y^{(1)} \in Y} f(x,y^{(0)}) f(x,y^{(1)})|^{1/2}.\end{equation} We apply this three times. At the first application we take $X := \{y_2, y_3\}$ and $Y = \{y_1\}$, and put the function $f_1$ inside the $\b()$ term. We now have variables $y_1^{(0)}, y^{(1)}_1, y_2, y_3$. Now set $X := \{y_1^{(0)}, y^{(1)}_1, y_3\}$, $Y := \{y_2\}$ and arrange for everything involving $f_2$ to be placed in the $\b()$ term. We now have variables $y_1^{(0)}, y^{(1)}_1, y^{(0)}_2, y^{(1)}_2, y_3$. For the final application of Cauchy-Schwarz set $X := \{y^{(0)}_1, y^{(1)}_1, y^{(0)}_2, y^{(1)}_2\}$ and $Y := \{y_3\}$, and arrange for everything involving $f_3$ to be placed in the $\b()$ term. Note that at this point we have eliminated everything involving $f_1, f_2, f_3$ and have
\begin{align*} |&\E_{y_1, y_2, y_3} \textstyle f_1(-\frac{1}{2}  y_2 - 2y_3) f_2(\textstyle\frac{1}{3} y_1 - y_3) f_3(\frac{2}{3} y_1 + \frac{1}{2}y_2) f_4(y_1 + y_2 + y_3)| \\ & \leq |\E_{y_1^{(0)}, y_1^{(1)}, y_2^{(0)}, y_2^{(1)}, y_3^{(0)}, y_3^{(1)}} f_4(y^{(0)}_1 + y^{(0)}_2 + y^{(0)}_3) f_4(y^{(1)}_1 + y^{(0)}_2 + y^{(0)}_3) \times \dots \\ & \qquad\qquad\qquad \dots \times f_4(y^{(1)}_1 + y^{(1)}_2 + y^{(1)}_3)|^{1/8}.\end{align*}
The right-hand side here is precisely $\Vert f_4 \Vert_{U^3}$. 

To show that $\Lambda(f_1,f_2,f_3,f_4)$ is bounded by the other expressions $\Vert f_i \Vert_{U^3}$, one may proceed similarly. We leave the details to the reader.\endproof

We now come to the central question of quadratic Fourier analysis: when is $\Vert f \Vert_{U^3}$ large? The first key observation is that the answer is not simply the same as for the $U^2$-norm.

\begin{lemma}[Key example]\label{key-ex}
There is a function $f : G \rightarrow \mathbb{C}$ with $\Vert f \Vert_{\infty} \leq 1$ such that $\Vert f \Vert_{U^3} = 1$, but such that $\Vert \widehat{f} \Vert_{\infty} \leq N^{-1/2}$.
\end{lemma}
\proof Before embarking on the proof, we must remark that $\Vert \cdot \Vert_{U^3}$ has only been defined for real-valued functions thus far. To define it for complex-valued functions, one must take complex conjugates of the terms $f(x+h_1)$, $f(x+h_2)$, $f(x+h_3)$ and $f(x+h_1 + h_2 + h_3)$. The extension to complex-valued functions facilitates the discussion of examples, but is not otherwise essential in the theory. Keeping track of complex conjugates is rather a tedious affair, so will endeavour to work with real functions whenever possible.

Set $f(x) = \omega^{x^T x}$.  We have
\[ \Vert f \Vert_{U^3}^8 = \E_{x,h_1,h_2,h_3} \omega^{x^T x - (x + h_1)^T(x+h_1) - \dots - (x+h_1 + h_2 + h_3)^T(x + h_1 + h_2 + h_3)} = 1.\] This can be seen by intelligent direct computation (or even by na\"{\i}ve direct computation); the phase vanishes since it is essentially the third derivative of a quadratic.

To evaluate $\Vert \widehat{f} \Vert_{\infty}$, observe that we have
\[ |\E_{x \in G} \omega^{x^Tx + r^Tx}| = |\prod_{j = 1}^n \E_{x_j \in \F_5} \omega^{x_j^2 + r_j x_j}| = 5^{-n/2}.\]
This concludes the proof of the lemma.\endproof

We conclude this first lecture by stating three key results in additive combinatorics which we will need in the second lecture. These results will all be discussed and proved in other lectures in this school. In these results, $0 < c < 1 < C$ are absolute constants.

\begin{proposition}[The Balog-Szemer\'edi-Gowers theorem]\label{bsg}
Let $G$ be an abelian group, and suppose that $A \subseteq G$ is a set with $|A| = n$. Suppose that there are at least $\delta n^3$ \emph{additive quadruples} in $A$, that is to say solutions to $a_1 + a_2 = a_3 + a_4$. Then there is a subset $A' \subseteq A$ with $|A' | \geq c\delta^{C} |A|$ such that $|A' + A'| \leq C\delta^{-C}|A'|$.
\end{proposition}
This result will be the subject of Antal Balog's lecture at the school.

\begin{proposition}[Freiman's theorem in finite fields]\label{frei-fin}
Let $p$ be a prime, and write $\F_p^n$ for the $n$-dimensional vector space over the finite field with $p$ elements. Suppose that $A \subseteq \F_p^n$ is a set with $|A + A| \leq K|A|$. Then there is a subspace $H \leq \F_p^n$ such that $A \subseteq H$ and for which we have the bound $|H| \leq p^{CK^C}|A|$. 
\end{proposition}
This result will be discussed by Imre Ruzsa.

\scriptsize

\noindent\textit{Exercises.} For the reader wishing to familiarise herself with the Gowers norms, we offer a handful of exercises. Discussions pertinent to these exercises may be found in the papers \cite{gowers-longaps,green-tao-longprimeaps,green-tao-u3inverse}.

\textbf{1.} Let $k \geq 2$ be any integer, and define the Gowers $U^k$-norm by
\begin{equation}\label{us-def} \Vert f \Vert_{U^k}^{2^k} := \E_{x,h_1,\dots,h_k \in G} \prod_{\omega \in \{0,1\}^k} f(x + \omega \cdot h).\end{equation} Show that $\Vert \cdot \Vert_{U^k}$ is a norm. (\emph{Hint:} first define the \emph{Gowers inner product} $\langle f_{\omega} \rangle_{\omega \in \{0,1\}^k}$ for $2^k$ functions $(f_{\omega})_{\omega \in \{0,1\}^k}$ by modifying \eqref{us-def}. Then use several applications of the Cauchy-Schwarz inequality to prove the \emph{Gowers-Cauchy-Schwarz inequality} 
\[ |\langle f_{\omega} \rangle_{\omega \in \{0,1\}^k}| \leq \prod_{\omega} \Vert f_{\omega} \Vert_{U^k}.\]
Finally, use this to prove the triangle inequality for $\Vert f \Vert_{U^k}$).

\textbf{2.} Prove that the Gowers $U^k$-norms are \emph{nested}:
\[ \Vert f \Vert_{U^2} \leq \Vert f \Vert_{U^3} \leq \dots.\]

\textbf{3.} By generalising Lemma \ref{key-ex}, show that the Gowers norms are \emph{strictly} nested in the following strong sense. For any $k \geq 3$ there is $c_k > 0$ such that the following is true. For any $N$, there is a group $G$ with $|G| \geq N$ and a function $f : G \rightarrow \mathbb{C}$ with $\Vert f \Vert_{\infty} = \Vert f \Vert_{U^k} = 1$ such that $\Vert f \Vert_{U^{k-1}} \ll N^{-c_k}$.

\textbf{4.} We noted that the $U^2$ inverse theorem is an if and only if statement. That is, if $f$ is a bounded function with $|\E_x f(x)\omega^{r^T x}| \geq \delta$ for some $r$ then $f$ has large $U^2$-norm. Prove this without using the fact that $\Vert f \Vert_{U^2} = \Vert \widehat{f} \Vert_4$. (\emph{Hint}: use the Gowers-Cauchy-Schwarz inequality of Ex. 1.)

\textbf{5.} Let $G = \F_5^n$. Suppose that 
\[ | \E_{x \in G} f(x) \omega^{x^T M x + r^t x} | \geq \delta.\] for some matrix $M$ and vector $r$. Prove that $\Vert f \Vert_{U^3} \geq \delta$. (\emph{Hint}: apply the Gowers-Cauchy-Schwarz inequality again. You will need the generalisation of $\Vert \cdot \Vert_{U^3}$ which covers complex-valued functions; this can be obtained by inserting appropriate complex conjugate symbols, as was discussed during the proof of Lemma \ref{key-ex}.

\textbf{6.} (Generalising the generalised von Neumann theorem) Show that
\[ |\Lambda_k(f_1,\dots,f_k)| \leq \inf_{i=1,\dots,k} \Vert f_i \Vert_{U^{k-2}}.\]

\noindent\textit{Further reading.} This material was originally laid out in Gowers \cite{gowers-longaps}, though the notation was slightly different and (of course) the Gowers norms were not named as such! Various expositions of the material may be found in papers by one or both of Terry Tao and myself. See, for example, \cite{green-tao-longprimeaps,green-tao-u3inverse}.

A very general version of the generalised von Neumann theorem (linking systems of $s$ equations in $t$ unknowns to the $U^{s+1}$ norm) may be found in our forthcoming paper \cite{green-tao-u3mobius}, and an even more general version (applying to functions which are not necessarily bounded by 1) may be found in \cite{green-tao-linearprimes}.

Analogues of much of the material in this lecture were discovered in ergodic theory about 20 years ago. For more on this fascinating connection, the lectures of Kra in these Proceedings are illuminating. 

The Balog-Szemer\'edi-Gowers theorem was originally proved by Gowers \cite{gowers-4aps}, and is a quantitative version of the earlier result of Balog and Szemer\'edi \cite{balog} (see also Balog's article in these Proceedings). A version with a good value of the exponent $C$ may be found in \cite{chang-er}. This material is also covered in my notes \cite{green-bkt}. The Pl\"unnecke--Ruzsa inequality was obtained in \cite{plun} and afforded an elegant proof by Ruzsa in \cite{ruzsa-graph}. The original reference for Proposition \ref{frei-fin} is the paper \cite{ruzsa-frei} by Imre Ruzsa. For self-contained notes on Pl\"unnecke's inequality and Freiman's theorem, see \cite{green-mit}. For a discussion of all of the material in this lecture (and indeed much of the material in the other lectures) see the book \cite{tao-vu-book}.

\normalsize

\section{Lecture 2}

Topics to be covered:

\begin{itemize}
\item The inverse theorem for the $U^3$-norm on $\F_5^n$.
\end{itemize}

\emph{Some notation.} Let $E,E'$ be real-valued expressions. We will write $E \gg_{\delta} E'$ to mean that there is some function $c(\delta) > 0$ such that $E \geq c(\delta) E'$. There is nothing particularly unusual about this notation, but one aspect of the manner in which we shall apply it is somewhat subtle. When we write, for example, ``let $N \gg_{\delta} 1$'', we mean ``let $N \geq c(\delta)$, where $c : \R_+ \rightarrow \R_+$ is some function which may be chosen so that later arguments work''. We do \emph{not} (of course) mean that an arbitrary function $c$ may be chosen.

We will also, on occasion, use the notation $O_{\delta}(1)$ to denote a finite quantity which depends only on $\delta$.

We have deliberately chosen topics within the subject of quadratic Fourier analysis for which bounds are unimportant, since these are the topics most allied to the ``infinitary'' ideas which feature in the lectures of Kra and Tao in these Proceedings. It is quite reasonable to think of there being just two types of quantity in these lectures: \emph{finite} quantities which depend only on $\delta$, and \emph{infinite} quantities which depend on the size of $\F_5^n$.

Let us recall the main question we are trying to address.

\begin{question}[Gowers inverse question] Suppose that $f : G \rightarrow [-1,1]$ is a function and that $\Vert f \Vert_{U^3} \geq \delta$. What can we say about $f$?
\end{question}

It turns out to be \emph{much} easier to address this question in a finite field setting such as $G = \F_5^n$. We showed in the exercises to Lecture 1 that if $f$ correlates with a quadratic phase $\omega^{x^T M x + r^T x}$ then $f$ has large $U^3$ norm. It turns out that the converse is also true, though this is much harder to prove and will be our main goal in this lecture.

\begin{proposition}[Inverse theorem for the $U^3$-norm on $\F_5^n$]\label{inv-thm-1} Suppose that $f : G \rightarrow [-1,1]$ is a function for which $\Vert f \Vert_{U^3} \geq \delta$. Then there is a matrix $M \in \GL_n(\F_5)$ and a vector $r \in \F_5^n$ so that
\[ |\E_{x \in G} f(x)\omega^{x^T M x + r^T x}| \gg_{\delta} 1.\]
\end{proposition}
\emph{Remark.} Write $E := \sup_{r,M} |\E_{x \in G}f(x) \omega^{x^T M x + r^T x}|$. It is not hard to check that the proof we give would allow one to replace $E \gg_{\delta} 1$ by some bound of the form $E \geq \exp(-C\delta^{-C})$. For our later application, we will merely need \emph{some} lower bound of the form $E \gg_{\delta} 1$. There are other applications where bounds are important -- see the \emph{further reading} at the end of this lecture for a discussion.

To prove Proposition \ref{inv-thm-1} we will essentially follow the approach of Gowers \cite{gowers-4aps}. We will, however, employ a slight twist which is essentially due to Samorodnitsky \cite{samorod}. 

\begin{definition}[Derivatives]
Suppose that $f : G \rightarrow \R$ is a function. Then for any $h \in G$ we define the function $\Delta(f;h)$ by
\[ \Delta(f;h)(x) := f(x) f(x-h).\]
\end{definition}

\emph{Remark.} It is convenient, though perhaps slightly mystifying, to give the name ``derivative'' to this construction. If we extended the definition to complex-valued functions by setting $\Delta(f;h)(x) = f(x)\overline{f(x-h)}$ and applied it with $f(x) = e^{2\pi i \phi(x)}$, the mystery might be reduced somewhat as the phase $\phi$ is indeed being differentiated.

\begin{proposition}[Samorodnitsky's identity]
Let $f : G \rightarrow \R$ be any function. Then we have
\begin{equation}\label{eq40} \sum_{r_1 + r_2 = r_3 + r_4} \E_{h_1 + h_2 = h_3 + h_4} |\Delta(f;h_1)^{\wedge}(r_1)|^2 \dots |\Delta(f;h_4)^{\wedge}(r_4)|^2 = \E_h\Vert \Delta(f;h)^{\wedge}\Vert_8^8.\end{equation}
\end{proposition}
\proof The idea of the proof is simple: we show that both sides are equal to 
\begin{equation}\label{star} \sum_{(c_1,\dots,c_8,c'_1,\dots,c'_8) \in \mathcal{C}} f(c_1) \dots f(c_8) f(c'_1) \dots f(c'_8),\end{equation}
where the sum is over all configurations $\mathcal{C}$ with
\[ c_1 + \dots + c_4 = c_5 + \dots + c_8\] and \[ c'_1 - c_1 = \dots = c'_8 - c_8.\]
To show that the RHS of \eqref{eq40} is equal to \eqref{star} is the easier of the two tasks to accomplish. One notes that
\[ \|\Delta(f;h)^{\wedge}\Vert_8^8 = \E_x |\Delta(f;h) \ast \Delta(f;h) \ast \Delta(f;x) \ast \Delta(f;h)(x)|^2,\] by Parseval's identity and the fact that $(f \ast g)^{\wedge} = \widehat{f}\widehat{g}$. That the expectation of this over $h$ is equal to \eqref{star} follows by expansion.

To prove that the LHS of \eqref{eq40} is equal to \eqref{star}, it is convenient to introduce some notation. If $\psi : \widehat{G} \rightarrow \C$ is a function then we define $\psi^{\vee} : G \rightarrow \C$ by
\[ \psi^{\vee}(x) := \sum_{r \in \widehat{G}} \psi(r) \omega^{-r^T x}.\]
Note that the inversion formula is equivalent to 
\begin{equation}\label{inv-new} (\widehat{f})^{\vee} = f.\end{equation}
If $\psi,\phi : \widehat{G} \rightarrow \C$ are two functions then we define
\[ \psi \ast \phi(r) := \sum_{s \in \widehat{G}} \psi(s) \phi(r - s)\] and note the formula
\[ (\psi \ast \phi)^{\vee} = \psi^{\vee}\phi^{\vee}.\]
It follows from these facts and Parseval's identity that for any four functions $g_1,\dots,g_4 : G \rightarrow \C$ we have
\begin{equation}\label{star-88} \sum_{r_1 + r_2 = r_3 + r_4} \widehat{g_1}(r_1)\widehat{g_2}(r_2)\overline{\widehat{g_3}(r_3)\widehat{g_4}(r_4)} = \sum_{r} \widehat{g_1} \ast \widehat{g_2}(r) \overline{\widehat{g_3} \ast \widehat{g_4}(r)} = \E_x g_1(x)g_2(x) \overline{g_3(x)g_4(x)}.\end{equation}
We apply this with
\[ g_i = \Delta(f;h_i) \ast \Delta(f;h_i)^{\circ},\]
where we have defined $f^{\circ}(x) := \overline{f(-x)}$. Noting that $(f^{\circ})^{\wedge} = \overline{\widehat{f}}$, we see that 
\[ \widehat{g_i}(r) = |\Delta(f;h_i)^{\wedge}(r)|^2.\]
Substituting into \eqref{star-88}, we see that the LHS of \eqref{eq40} is equal to
\[ \E_{h_1 + h_2 = h_3 + h_4} \E_x \prod_{i=1}^4 \Delta(f;h_i) \ast \Delta(f;h_i)^{\circ}(x).\] Expanding out, we recover \eqref{star} once more.\endproof

Using this identity, we can prove the following crucial result, which provides the first link between functions $f$ with large $U^3$-norm and quadratic phases. It states that the derivatives $\Delta(f;h)$ obey a sort of weak linearity property.

\begin{proposition}[Gowers]\label{gowers-prop}  Let $f : G \rightarrow [-1,1]$ be a function, and suppose that $\Vert f \Vert_{U^3} \geq \delta$. Suppose that $|G| \gg_{\delta} 1$. Then there is a function $\phi : G \rightarrow \widehat{G}$ such that 
\begin{enumerate}
\item $|\Delta(f;h)^{\wedge}(\phi(h))| \gg_{\delta} 1$ for all $h \in S$, where $|S| \gg_{\delta} |G|$;
\item There are $\gg_{\delta} |G|^3$ quadruples $(s_1,s_2,s_3,s_4) \in S^4$ such that $s_1 + s_2 = s_3 + s_4$ and $\phi(s_1)+\phi(s_2) = \phi(s_3)+\phi(s_4)$.
\end{enumerate}
\end{proposition}
\proof Set $N := |G|$. One may easily check that
\[ \Vert f \Vert_{U^3}^8 = \E_h \Vert \Delta(f;h)\Vert_{U^2}^4.\] Recalling that the $U^2$-norm is the $L^4$ norm of the Fourier transform, we thus have
\[ \Vert f \Vert_{U^3}^8 = \E_h \|\Delta(f;h)^{\wedge}\|_4^4.\]
Now H\"older's inequality and Parseval's identity imply that for any $h$ we have\[ \Vert\Delta(f;h)^{\wedge}\Vert_4^4 \leq \Vert\Delta(f;h)^{\wedge}\Vert_2^{4/3} \Vert \Delta(f;h)^{\wedge}\Vert_8^{8/3} \leq \Vert \Delta(f;h)^{\wedge}\Vert_8^{8/3}.\]
Another application of H\"older yields
\[ \E_h  \| \Delta(f;h)^{\wedge}\|^{8/3}_8 \leq \big( \E_h \| \Delta(f;h)^{\wedge}\|_8^8  \big)^{1/3}.\]
Combining these observations, we conclude that
\[ \E_h  \| \Delta(f;h)^{\wedge}\|^8_8  \geq \delta^{24}.\] Samorodnitsky's identity then allows us to conclude that
\begin{equation}\label{eq778} \sum_{r_1 + r_2 = r_3 + r_4} \E_{h_1 + h_2 = h_3 + h_4} |\Delta(f;h_1)^{\wedge}(r_1)|^2 \dots |\Delta(f;h_4)^{\wedge}(r_4)|^2 \geq \delta^{24} .\end{equation}
To each $h \in G$, we associate the set $\Phi(h)$ of characters $r$ for which $|\Delta(f;h)^{\wedge}(r)| \geq \delta^{50}$. It is immediate from Parseval's identity that $|\Phi(h)| \leq \delta^{-100}$ for all $h$. Now the contribution to \eqref{eq778} from those $h_i,r_i$ for which $r_1 \notin \Phi(h_1)$ (say) is bounded by
\[ \delta^{100}\sum_{r_2,r_3,r_4} \E_{h_2,h_3,h_4} |\Delta(f;h_2)^{\wedge}(r_2)|^2 |\Delta(f;h_3)^{\wedge}(r_3)|^2 |\Delta(f;h_4)^{\wedge}(r_4)|^2 \leq \delta^{100}. \]
It follows that
\[ \sum_{r_1 + r_2 = r_3 + r_4} \E_{h_1 + h_2 = h_3 + h_4} 1_{r_1 \in \Phi(h_1)}|\Delta(f;h_1)^{\wedge}(r_1)|^2 \dots 1_{r_4 \in \Phi(h_4)}|\Delta(f;h_4)^{\wedge}(r_4)|^2 \geq \delta^{24}/2,\] and so in particular there are at least $\delta^{24} N^3/2$ \emph{additive octuples} $(h_1,r_1,\dots,h_4,r_4)$ such that $h_1 + h_2 = h_3 + h_4$, $r_1 + r_2 = r_3 + r_4$ and $r_i \in \Phi(h_i)$ for $i = 1,\dots,4$. We say that an octuple is \emph{proper} if $h_1,\dots,h_4$ are all distinct. The number of our additive octuples which fail to be proper is clearly $\ll_{\delta} N^2$ and hence, since $N$ is so large, at least $\delta^{24}N^3/4$ of them \emph{are} proper.

Let $S$ be the set of all $h$ for which $\Phi(h) \neq \emptyset$. It is easy to see that $|S| \gg_{\delta} |G|$, since otherwise there could not be enough additive octuples. For each $h \in S$, pick an element $\phi(h)$ uniformly at random from $\Phi(h)$, and suppose that these choices are independent for different $h$. For each proper additive octuple $(h_1,r_1,\dots,h_4,r_4)$, the probability that  it \emph{fits $\phi$}, that is to say that $r_i = \phi(h_i)$ for $i= 1,2,3,4$, is precisely $1/|\Phi(h_1)| \dots |\Phi(h_4)|$. This is $\gg_{\delta} 1$. It follows that the expected number of additive octuples which fit $\phi$ is $\gg_{\delta} |G|^3$. In particular there is some specific choice of $\phi$ for which $\gg |G|^3$ additive octuples fit $\phi$.

It takes a few seconds to realise that we have, in fact, proved the result. Indeed, an octuple which fits $\phi$ is precisely an additive quadruple of points $h_1,\dots,h_4$ such that $\phi(h_1)+\phi(h_2) = \phi(h_3)+\phi(h_4)$ and $\phi(h_i) \in \Phi(h_i)$, that is to say $|\Delta(f,h)^{\wedge}(\phi(h))| \geq \delta^{50}$.\endproof 

We have made a crucial step: assuming that $\Vert f \Vert_{U^3}$ was large, we deduced that the derivative of $f$ has a certain weak linearity property. We must now work with this property and make it somewhat stronger.

\begin{proposition}[From weak linearity to linearity]\label{wk-to-strong}
Suppose that $\phi : G \rightarrow \widehat{G}$ is a function with the property in Proposition \ref{gowers-prop} \textup{(2)}, that is to say there is some set $S \subseteq G$ with $|S| \gg_{\delta} |G|$ such that there are $\gg_{\delta} |G|^3$ additive quadruples $(s_1,s_2,s_3,s_4)$ such that $s_1 + s_2 = s_3 + s_4$ and $\phi(s_1) + \phi(s_2) = \phi(s_3) + \phi(s_4)$. Then there is some linear function $\psi(x) = Mx + b$, where $M \in \GL_n(\F_5)$ and $b \in \F_5^n$, such that $\phi(x) = \psi(x)$ for $\gg_{\delta} |G|$ values of $x \in S$.
\end{proposition}
\proof The first step is to observe that the conclusion of Proposition \ref{gowers-prop} may be rephrased using the \emph{graph}
\[ \Gamma := \{ (h,\phi(h)) : h \in S\},\] which is a subset of $G \times \widehat{G}$. Statement (2) of Proposition \ref{gowers-prop} is just the same as saying that $\Gamma$ has $\gg_{\delta} |G|^3$ additive quadruples. It follows from the Balog-Szemer\'edi-Gowers theorem that there is a subset $\Gamma' \subseteq \Gamma$ with \[ |\Gamma'| \gg_{\delta}|\Gamma| \gg_{\delta} |G|\]  and \[ |\Gamma' + \Gamma'| \ll_{\delta} |\Gamma'|.\]
Define $S' \subseteq S$ by
\[ \Gamma' := \{(h,\phi(h)) : h \in S'\},\]
and note that 
\[ |S'| \gg_{\delta} |G|.\]
Now we may identify $G \times \widehat{G}$ with $\F_5^n \times \F_5^n$ and hence with $\F_5^{2n}$. From Ruzsa's finite field analogue of Freiman's theorem, it follows that there is some subspace $H \leq \F_5^n \times \F_5^n$, 
\begin{equation}\label{h-upper} |H| \ll_{\delta} |G|,\end{equation} such that $\Gamma' \subseteq H$.

Consider the map $\pi : H \rightarrow G$ onto the first factor. The image of this linear map contains $S'$, and so from \eqref{h-upper} and the lower bound for $|S'|$ we see that
\[ \dim_{\F_5}\ker \pi \ll_{\delta} 1.\]
It follows that we may foliate $H$ into $\ll_{\delta} 1$ cosets of some subspace $H'$, such that $\pi$ is injective on each of these cosets. By averaging, we see that there is some $x$ such that
\[ |(x + H') \cap \Gamma'| \gg_{\delta} |G|.\]
Set $\Gamma'' := (x + H') \cap \Gamma'$, and define $S'' \subseteq S'$ accordingly. Then $\pi|_{x + H'}$ is an affine isomorphism onto its image $V$, which means that there is an affine linear map $\psi : V \rightarrow \widehat{G}$ such that $(s'',\psi(s'')) \in \Gamma''$ for all $s'' \in S''$, that is to say  $\psi(s'') = \phi(s'')$ for all $s'' \in S''$.\endproof

Let us put this last result together with Proposition \ref{gowers-prop}.

\begin{corollary}[Linearity of the derivative]\label{lin-deriv-cor} Suppose that $f : G \rightarrow [-1,1]$ is a function with $\Vert f \Vert_{U^3} \geq \delta$. Suppose that $|G| \gg_{\delta} 1$. Then there is some $M \in \GL_n(\F_5)$ and some $b \in \F_5^n$ such that 
\[ \E_h |\Delta(f;h)^{\wedge}(M h + b)|^2 \gg_{\delta} 1.\]
\end{corollary}
\proof Recall that $\phi$ is defined for $h \in S$, where 
\[ |S| \gg_{\delta} |G|\] and that it has the property that
\[ |\Delta(f;h)^{\wedge}(\phi(h))| \gg_{\delta} 1\] 
for all $h \in S$. We proved in Proposition \ref{wk-to-strong} that there is an affine linear function $\psi(h) = Mh + b$ such that $\phi(h) = \psi(h)$ for all $h \in S''$, where $|S''| \gg_{\delta} |G|$. The corollary follows immediately.\endproof

Corollary \ref{lin-deriv-cor} shows that the derivative of a function $f$ with large $U^3$ norm correlates with a linear function. Recall that our aim is to show that $f$ correlates with a quadratic function $x \mapsto \omega^{x^T M x + r^T x}$. This latter function does have a linear derivative, but this derivative is \emph{symmetric}. For that reason we need the following lemma, which states that the matrix $M$ in Corollary \ref{lin-deriv-cor} is automatically nearly symmetric.

\begin{lemma}[Symmetry argument]\label{sym-lem}
Suppose that $f : G \rightarrow [-1,1]$ is a function, that $M \in \GL_n(\F_5)$,and that $b \in \F_5^n$. Suppose that 
\[ \E_h |\Delta(f;h)^{\wedge}(Mh + b)|^2 \gg_{\delta} 1.\] Then $M$ is approximately symmetric in the sense that
\[  \rk(M - M^T) \ll_{\delta} 1.\]
\end{lemma}
\proof  Write $D = M - M^{T}$. Expanding the assumption gives
\[ \E_{x,y,h} f(x)f(x-h)f(y)f(y-h)\omega^{(x-y)^TMh + (x-y)^Tb} \gg_{\delta} 1,\] Making the substitution $z = x + y - h$, this becomes
\[ \E_{x,y,z} f(x)f(z-x)f(y)f(z-y)\omega^{(x-y)^TM(x + y - z) + (x-y)^Tb} \gg_{\delta} 1,\] which can be written
\[ \E_z \E_x \Delta'(f;z)(x)\omega^{x^TM(x-z) + x^Tb}\E_y \Delta'(f;z)(y)\omega^{-y^TM(y - z) - y^Tb}\omega^{x^TDy} \gg_{\delta} 1.\]
Here, we have written \[ \Delta'(f;z)(t) := f(t)f(z - t).\] Writing 
\[ g_z(x) := \Delta'(f;z)(x)\omega^{x^TM(x-z) + x^Tb},\]
we have
\[ \E_z \E_{x,y} g_z(x)\overline{g_z(y)}\omega^{x^TDy} \gg_{\delta} 1.\]
Averaging over $z$, we see that there is some function $g : G \rightarrow \C$ with $\Vert g \Vert_{\infty} \leq 1$ such that
\[ |\E_x g(x)\overline{g(y)} \omega^{x^T D y}| \gg_{\delta} 1,\] that is to say
\[ |\E_x g(x) \widehat{g}(Dx)| \gg_{\delta} 1.\]
This implies that
\[ \E_x |\widehat{g}(Dx)| \gg_{\delta} 1, \] and so in particular there are $\gg_{\delta} |G|$ values of $x$ such that $|\widehat{g}(Dx)| \gg_{\delta} 1$. However we know from Parseval's identity that the number of $r$ such that $|\widehat{g}(r)| \gg_{\delta} 1$ is $\ll_{\delta} 1$. Thus there is some set $S \subseteq \F_5^n$ with $|S| \gg_{\delta} |G|$ and $|D(S)| \ll_{\delta} 1$. This implies that 
\[ |\ker(D)| \gg_{\delta} |G|,\] which immediately implies the result.\endproof

We have shown that if $\Vert f \Vert_{U^3}$ is large then the derivative of $f$ correlates with a symmetric linear form. To complete the proof of Proposition \ref{inv-thm-1}, we must ``integrate'' this statement and show that $f$ correlates with a quadratic. We give this integration now.

\emph{Proof of Proposition \ref{inv-thm-1}.} From Corollary \ref{lin-deriv-cor} and Lemma \ref{sym-lem}, we know that 
\begin{equation}\label{assump-3} \E_h |\Delta(f;h)^{\wedge}(Mh + b)|^2 \gg_{\delta} 1,\end{equation}
where 
\[ \rk(M - M^T) \ll_{\delta} 1.\]
Write $M_{\sym} := \frac{1}{2}(M + M^T)$, and let $V := \ker(M - M^T)$. For each $t \in G$ there is some $b_t$ such that we have
\[ Mh + b = M_{\sym} h + b_t.\]
for all $h \in V + t$. By a trivial averaging argument and the fact that $\codim(V) \ll_{\delta} 1$, we may find a $t$ such that
\[ \E_h 1_{h \in V + t}|\Delta(f;h)^{\wedge}(Mh + b)|^2 \gg_{\delta} 1.\]
This of course implies that
\[ \E_h 1_{h \in V + t} |\Delta(f;h)^{\wedge}(M_{\sym}h + b_t)|^2 \gg_{\delta} 1,\]
and hence by positivity that
\[ \E_h |\Delta(f;h)^{\wedge}(M_{\sym} h + b_t)|^2 \gg_{\delta} 1.\]
By redefining $M$ to be $M_{\sym}$ and $b$ to be $b_t$, it follows that we may assume in \eqref{assump-3} that $M$ is symmetric.

Expanding out \eqref{assump-3} we obtain
\[ \E_{h,x,y} f(x)f(x-h)f(y)f(y-h)\omega^{h^T M (x-y) + b^T (x-y)} \gg_{\delta} 1.\]
Substituting $y := x - k$, we obtain
\[ \E_{h,x,k} f(x) f(x-h)f(x-k)f(x-h-k) \omega^{h^T M k + b^T k} \gg_{\delta} 1.\]
Using the identity
\[ x^T M x - (x - h)^T M (x-h) - (x - k)^T M (x-k) + (x - h - k)^TM(x - h - k) = 2h^T M k,\]
 this may be written as
\begin{equation}\label{nearly-there} \E_{h,x,k} g_1(x) g_2(x-h) g_3(x - k)g_4(x - h - k) \gg_{\delta} 1,\end{equation}
where $g_1(x) := f(x) \omega^{\frac{1}{2}x^T M x}$, $g_2(x) := f(x)\omega^{-\frac{1}{2}x^T M x - b^Tx}$, $g_3(x) :=  f(x)\omega^{-\frac{1}{2}x^T M x}$ and  $g_4(x) := f(x)\omega^{\frac{1}{2}x^T M x - b^Tx}$. Note that the functions $g_2,g_3,g_4$ are bounded by 1; this is, in fact, the only property of them that we shall use.

Now the left-hand side of \eqref{nearly-there} may be rewritten using the Fourier transform as
\[ \sum_r \widehat{g_1}(r) \widehat{g_2}(-r) \widehat{g_3}(-r) \widehat{g_4}(r).\]
It follows immediately frrm H\"older's inequality that \[
\Vert \widehat{g}_1 \Vert_4 \gg_{\delta} 1,\]
which, since $\Vert \widehat{g}_1 \Vert_2 \leq 1$, implies that
\[
\Vert \widehat{g}_1 \Vert_{\infty} \gg_{\delta} 1,\]that is to say there is some $r \in \F_5^n$ such that
\[ |\E_x f(x) \omega^{\frac{1}{2} x^T M x + r^T x} | \gg_{\delta} 1.\]
This, at last, completes the proof of Proposition \ref{inv-thm-1}.\endproof

\emph{Remark.} In going from \eqref{nearly-there} to the end of the proof, what we have really done is apply the Gowers-Cauchy-Schwarz inequality (cf. the exercises following Lecture 1) and the inverse theorem for the $U^2$-norm.

\scriptsize
\emph{Further reading.} The orginal argument of Gowers is in \cite{gowers-4aps}. This took place in the group $G = \Z/N\Z$, not in a finite field model, and did not quite give a necessary and sufficient inverse theorem for the $U^3$-norm. It was instead shown that if $f : \Z/N\Z \rightarrow [-1,1]$ has large $U^3$-norm then $f$ correlates with a quadratic polynomial on some subprogression of length a power of $N$. This is a ``local'' statement, and as such is much weaker than having large $U^3$-norm, which is ``global'', i.e. involves averaging over the whole group $G$.

To get an inverse theorem, one extra ingredient must be added to Gowers' work. This is the symmetry argument, Lemma \ref{sym-lem}. It was first given in \cite{green-tao-u3inverse}. That paper gives an inverse theorem for the $U^3$-norm in any finite abelian group of odd order. To even state the result is somewhat complicated, and we defer a discussion until we have thoroughly examined the finite field case. An inverse theorem for the $U^3$-norm in $\F_2^n$ was given by Samorodnitsky \cite{samorod}, using the method we have described but with a slight twist to enable him to handle characteristic 2. It is very likely that a combination of his methods and ours would allow one to prove an inverse theorem in \emph{any} finite abelian $G$, but to my knowledge no-one has yet undertaken this task.

As we remarked, one may replace our $\gg_{\delta} 1$ notation with more precise bounds, ending up with a version of Proposition \ref{inv-thm-1} with a function of the form $\exp(-C\delta^{-C})$ on the right-hand side. It would be of great interest to know whether this could be improved, perhaps even to $c\delta^C$. This would follow from the so-called Polynomial-Freiman-Ruzsa conjecture, the finite field version of which is discussed in \cite{green-fin-field}.

The strongest known inverse result for the $U^3$ norm on $\F_5^n$ is the following, proved in \cite{green-tao-u3inverse}.

\begin{proposition}[Inverse theorem for the $U^3$-norm on $\F_5^n$, II]\label{inv-thm-2} Suppose that $f : \F_5^n \rightarrow [-1,1]$ is a function for which $\Vert f \Vert_{U^3} \geq \delta$. Then there exists a subspace $H \leq \F_5^n$ with $\codim(H) \leq C\delta^{-C}$, together with a system of quadratic forms $r_y^T x + x^T M_y x$ indexed by the cosets $y + H$ of $H$, such that 
\[ \E_y|\E_{x \in y + H} f(x) \omega^{x^T M_y x + r_y^T x}| \geq c\delta^{C}.\]
\end{proposition}

Note that the amount of correlation is $c\delta^C$ rather than $\exp(-C\delta^{-C})$, but one must pass to a coset of a subspace of somewhat large codimension. 

The proof of this result is rather longer than that of Proposition \ref{inv-thm-1}, and involves a good deal more machinery (Bogolyubov's method and Freiman homomorphisms). This stronger result is necessary for certain applications, for example in our paper \cite{green-tao-ffszem} in which it is shown that $r_4(\F_5^n) \ll N(\log N)^{-c}$.
\normalsize

\section{Lecture 3}

Topics to be covered:

\begin{itemize}
\item Quadratic factors
\item The energy increment lemma
\item The idea of approximating a function by projecting onto a low-complexity factor
\item The Koopman-von Neumann decomposition
\item The arithmetic regularity decomposition
\end{itemize}

Our main effort so far has been devoted to proving a result of the form ``if $\Vert f \Vert_{U^3}$ is large then $f$ has a large quadratic Fourier coefficient''. 

In this section we turn to a discussion of how this kind of information can be useful to us. There are many instances in additive combinatorics where study of a single Fourier coefficient is fruitful. However there are many other occasions on which it is beneficial to consider \emph{several} Fourier coefficients of $f$, say the set of large Fourier coefficients of $f$.
We must develop analogues of this theory in the quadratic setting.

From now on, matrices $M \in \GL_n(\F_5)$ will only appear in quadratic forms $x^T M x$. Thus from this point onwards it is natural to adopt the convention that \emph{all matrices are symmetric}. We note that a (slightly) more high-brow approach to the whole theory, avoiding the use of bases, appears in our paper \cite{green-tao-ffszem}.

The following simple lemma will be used over and over again.

\begin{lemma}[Gauss sums]\label{gauss-sums}
Suppose that $M$ is symmetric and that $\rk M = d$. Then for any $r \in G$ we have
\[ |\E_{x \in G} \omega^{x^T M x + r^T x}| \leq 5^{-d/2}.\] If $r = 0$ then equality occurs.
\end{lemma}
\proof Squaring, we obtain
\begin{align*} |\E_{x \in G} \omega^{x^T M x + r^T x}|^2 &= \E_{h} \omega^{h^T M h + r^T h} \E_x \omega^{2h^T M x}\\ & \leq \E_h |\E_x \omega^{2h^T M x}|.\end{align*}
The inner sum is zero unless $h \in \ker(M)$. This occurs with probability $5^{-r}$, and so we do indeed get
\[ |\E_{x \in G} \omega^{x^T M x + r^T x}|^2 \leq 5^{-d}.\]
If $r = 0$ then the phase $\omega^{h^T M h + r^T h}$ is actually equal to 1 when $h \in \ker(M)$, and so equality occurs.
\endproof

Using this lemma, we may highlight one of the immediate difficulties with formulating ``quadratic Fourier analysis''.

\begin{lemma}[Profusion of large QFCs]
Let $f : \F_5^n \rightarrow [-1,1]$ be a function. Then there at most $\delta^{-2}$ values of $r$ for which 
\[ |\widehat{f}(r)| = |\E_{x \in \F_5^n} f(x) \omega^{r^T x}| \geq \delta.\] However, the number of pairs $(M,r)$ such that
\[ |\E_{x \in \F_5^n} f(x) \omega^{x^T M x + r^T x}| \geq \delta\]
need not be bounded in terms of $\delta$.
\end{lemma}
\proof The first statement, which is included for comparison with the classical setting, is immediate from Parseval's identity. To illustrate the second, one may consider a function as simple as $f(x) \equiv 1$. For any symmetric matrix $M$ with $\rk(M) \leq \log_5(1/\delta)$, we have
\[ |\E_{x \in \F_5^n} f(x) \omega^{x^T M x}| \geq \delta.\]
The number of such matrices is not bounded in terms of $\delta$.\endproof

This lemma suggests that we should perhaps only consider QFCs as ``essentially different'' if they are not too close in rank. This turns out to be a useful idea, and we will return to it later when we are in a position to formulate it properly.

As we said there are many arguments (e.g. \cite{croot-ap-lims,green-roth-primes,heath-brown-3aps,szemeredi-3aps}) where one considers the set of $\delta$-large Fourier coefficients \[ \Spec_{\delta}(f)  := \{r \in \F_5^n : |\widehat{f}(r)| \geq \delta\}.\] Without going into details of the applications, let us describe a useful way to think about the way this construction is often used.

\begin{definition}[Factors]
Let $\phi_1,\dots,\phi_k : \F_5^n \rightarrow \F_5$ be any functions. These functions describe a $\sigma$-algebra $\B$ on $\F_5^n$, the atoms of which are sets (of which there are at most $5^k$) of the form $\{x : \phi_1(x) = c_1,\dots, \phi_k(x) = c_k\}$. If $f : \F_5^n \rightarrow \C$ is a function then we often consider the \emph{conditional expectation} $\E(f | \B)$. Note that $\E(f|\B)(x)$ is just the average of $f$ over the atom $\B(x)$ which contains $x$. We will usually refer to $\sigma$-algebras arising in this way as \emph{factors}, by analogy with ergodic theory. We say that a factor $\B'$ \emph{refines} $\B$ if every atom of $\B'$ is contained in an atom of $\B$. Thus $\B'$ is at least as fine a partition of $\F_5^n$ as $\B$ is.
\end{definition}

\begin{definition}[Linear factors]
Suppose that $r_1,\dots,r_k \in \F_5^n$. Then the $\sigma$-algebra $\B$ whose atoms are the sets $\{x : r_i^T x = c_i, i = 1,\dots,k\}$ is called a \emph{linear factor} of complexity at most $k$.
\end{definition}

\begin{proposition}[Linear Koopman-von Neumann decomposition]\label{lin-kvn}
Let $f : \F_5^n \rightarrow [-1,1]$ be a function and let $\delta > 0$ be a parameter. Then there is a linear factor $\mathcal{B}$ of complexity at most $4\delta^{-4}$ such that 
\[ f = f_1 + f_2,\]
where 
\[ f_1 := \E(f | \mathcal{B})\]
and 
\[ \Vert f_2 \Vert_{U^2} \leq \delta.\]
\end{proposition}
\emph{Remark.} The Koopman-von Neumann theorem may be described in words as ``any bounded function is the sum of a ``low complexity'' function formed by projecting onto a linear factor, and a ``uniform'' function which is small in $U^2$. 

\proof The proof we give uses Fourier analysis, and does not generalise to give a result for the $U^3$-norm. We include it to justify the fact that this is a proposition which encodes the notion of ``taking all the large Fourier coefficients of $f$''.

Write $\eta := \delta^2/2$. Let $S := \Spec_{\eta}(f)$: note that by Parseval's identity we have $|S| \leq 4\delta^{-4}$. Let $H = S^{\perp}$ be the annihilator of $f$ and write $\mu_H$ for the Haar measure on $H$, that is to say $\mu_H := 1_H/\E1_H$. Define $f_1 := f \ast \mu_H$ and $f_2 := f - f \ast \mu_H$. It is not hard to see that $f_1 = \E(f | \B)$, where $\B$ is the factor defined by the linear functions $r^T x$, $r \in S$. To conclude the proof, we only need check that $\Vert \widehat{f}_2 \Vert_{\infty}$ is small. To that end, we have
\[ |\widehat{f}_2(r)| = |\widehat{f}(r)| |1 - \widehat{\mu}_H(r)|.\]
If $r \in \Spec_{\eta}(f)$ then $\widehat{\mu_H}(r) = 1$, and so $\widehat{f_2}(r) = 0$. If $r \notin \Spec_{\eta}(f)$ then by definition we have $|\widehat{f}(r)| \leq \eta$, and so $|\widehat{f}_2(r)| \leq 2\eta$ in this case. It follows that $\Vert \widehat{f_2} \Vert_{\infty} \leq 2\eta$, and thus by the inverse theorem for the $U^2$-norm we have
$\Vert f_2 \Vert_{U^2} \leq \sqrt{2\eta}$. The result follows.\endproof

\begin{definition}[Quadratic factors]\label{quad-fact-def}
Let $r_1,\dots,r_{d_1} \in \F_5^n$ be vectors, and let $M_1,\dots,M_{d_2}$ $\in \GL_n(\F_5)$ be symmetric matrices. We write $\B_1$ for the linear factor generated by the $r_j^T x$. Write $\B_2$ for the $\sigma$-algebra generated by the functions $r_j^T x$ and the pure quadratic functions $x^T M_j x$. Clearly $\B_2$ refines $\B_1$. We call the pair $(\B_1,\B_2)$ a (homogeneous) quadratic factor of complexity $(d_1,d_2)$.
\end{definition}

\begin{proposition}[Quadratic Koopman-von Neumann decomposition]\label{q-kvn} Let $(\B^{(0)}_1,\B^{(0)}_2)$ be a quadratic factor with complexity at most $(d_1^{(0)},d_2^{(0)})$.
Let $f : \F_5^n \rightarrow [-1,1]$ be a function and let $\delta > 0$ be a parameter. Then there is a quadratic factor $(\mathcal{B}_1,\mathcal{B}_2)$ of complexity at most $(d_1^{(0)} + O_{\delta}(1),d_2^{(0)} + O_{\delta}(1))$ which refines $(\B^{(0)}_1,\B^{(0)}_2)$, and such that 
\[ f = f_1 + f_2,\]
where 
\[ f_1 := \E(f | \mathcal{B}_2)\]
and 
\[ \Vert f_2 \Vert_{U^3} \leq \delta.\]
\end{proposition}

\emph{Remark.} For applications in which bounds are unimportant, it is better to apply the arithmetic regularity lemma which we will give later. A version of the Koopman-von Neumann theorem with reasonable bounds is the key tool in \cite{green-tao-ffszem}. In that application we take $(\B^{(0)}_1,\B^{(0)}_2)$ to be the trivial factor.

The key to proving the Koopman von Neumann decomposition lies in the following result.

\begin{lemma}[Energy increment]\label{energy-plus}
Let $(\B_1,\B_2)$ be a quadratic factor of complexity at most $(d_1,d_2)$, and let $f : \F_5^n \rightarrow [-1,1]$ be a function such that
\[ \| f - \E(f|\B_2)\|_{U^3} \geq \delta.\] Then exists a refinement $(\B'_1,\B'_2)$ of $(\B_1,\B_2)$
of complexity at most $(d_1+ 1, d_2+1)$ such that we have the \emph{energy increment}
\begin{equation}\label{en-inc} \| \E(f|\B'_2) \|_{2}^2 \geq \| \E(f|\B_2) \|_{2}^2 + c(\delta),\end{equation}
where $c : (0,1) \rightarrow \R_+$ is some non-decreasing function of $\delta$.
\end{lemma}

\proof The function $g := f - \E(f | \B_2)$ is certainly bounded by $2$, so we may apply the inverse theorem for the $U^3$-norm (Proposition \ref{inv-thm-1}) to conclude that there is a quadratic $x^T M x + r^T x$ so that
\begin{equation}\label{inv-thm-consequence} |\E_x g(x) \omega^{x^T M x + r^T x} | \geq c(\delta).\end{equation} We may clearly assume that $c : (0,1) \rightarrow \R_+$ is a non-decreasing function. The linear part $r^T x$ and the pure quadratic part $x^T M x$ of this quadratic together induce a quadratic factor $(\widetilde{\B}_1,\widetilde{\B}_2)$ of complexity $(1,1)$.

Now since $x^T M x + r^T x$ is $\widetilde{\B}_2$-measurable, it is clear that
\[ \E_x g(x) \omega^{x^T M x + r^T x} = \E_x \E(g|\widetilde{\B}_2)(x) \omega^{x^T M x + r^T x},\]
In particular, \eqref{inv-thm-consequence} implies that
\begin{equation}\label{itc2} \Vert \E(g | \widetilde{\B}_2)\Vert_1 \geq c(\delta).  \end{equation}
Now define $\B'_1 := \B_1 \vee \widetilde \B_1$ and $\B'_2 := \B_2 \vee \widetilde \B_2$. Again, the meaning of this is the obvious one; simply intersect all the atoms of $\B_i$ with those of $\widetilde \B_i$. It is clear that $(\B'_1,\B'_2)$ is a quadratic factor of complexity at most $(d_1 + 1,d_2 + 1)$. 

It remains to establish the energy increment \eqref{en-inc}. A key tool is

\noindent\textbf{Pythagoras' theorem.} Suppose that $\B,\B'$ are two $\sigma$-algebras on $\F_5^n$ such that $\B'$ refines $\B$. Let $f : \F_5^n \rightarrow [-1,1]$ be any function. Then
\[ \Vert \E(f | \B')\Vert_2^2 = \Vert \E(f | \B) \Vert_2^2 + \Vert \E(f | \B') - \E(f | \B)\Vert_2^2.\]

Now we have the chain of inequalities \begin{align*}
\| \E(f|\B'_2) \|_{2}^2 - \| \E(f|\B_2) \|_{2}^2
&= \| \E(f|\B'_2) - \E(f|\B_2) \|_{2}^2 \\
&= \| \E( g | \B'_2 ) \|_{2}^2 \\
&\geq \| \E( g | \widetilde \B_2 ) \|_{2}^2 \\
&\geq \| \E( g | \widetilde \B_2 ) \|_{1}^2 \\
&\geq c(\delta).
\end{align*}
The justification of these five lines uses respectively Pythagoras' theorem, the fact that $\B'_2$ refines $\B_2$, Pythagoras' theorem together with the fact that $\B'_2$ refines $\widetilde{\B}_2$, the Cauchy-Schwarz inequality, and \eqref{itc2}.\endproof

\emph{Proof of Proposition \ref{q-kvn}.} Start with $(\B_1,\B_2) = (\B_1^{(0)},\B^{(0)}_2)$. If
\begin{equation}\label{test} \Vert f - \E(f | \B_2) \Vert_{U^3} \leq \delta\end{equation} then \texttt{STOP}. 
Otherwise, we may apply Lemma \ref{energy-plus} to extend $(\B_1,\B_2)$ to a quadratic factor
with complexity incremented by at most $(1, 1)$ and the energy 
$\| \E(f|\B_2)\|_{2}^2$ incremented by at least $c(\delta)$. If \eqref{test} holds then \texttt{STOP}, otherwise repeat the process.
Since $f$ is bounded, the energy $\| \E(f|\B_2)\|_{2}^2$ lies in the interval $[0,1]$.  Since $c : (0,1) \rightarrow \R_+$ is non-decreasing, we cannot iterate the above 
procedure more than $1/c(\delta)$ times before we \texttt{STOP}.  The claim follows.\endproof

We will not give an application of the Koopman von-Neumann decomposition, since the interesting applications require quantitative versions of the result (cf. \cite{green-tao-ffszem}). The result has a significant shortcoming, which is that the uniformity parameter $\delta$ need not be small in terms of the complexity of $(\B_1,\B_2)$. For such situations there is another type of decomposition, which we call the \emph{arithmetic regularity lemma} because of an analogy with Szemer\'edi's regularity lemma in graph theory. We note that any use of this type of decomposition necessarily results in terrible ``tower-type'' bounds: see for example \cite{gowers-example,green-gafa}. As we have stated, however, bounds are not our concern in these lectures. 

\begin{proposition}[Arithmetic regularity lemma for $U^3$]\label{arith-reg-lem}
Let $\delta > 0$ be a parameter, and let $\omega : \R_+ \rightarrow \R_+$ be an arbitrary growth function\footnote{The use of arbitrary growth functions really does put us in the domain of ``discrete analogues of infinitary mathematics''. The arithmetic regularity lemma is indeed very close in spirit to the main result of the ergodic-theoretic paper \cite{bhk}.} \textup{(}which may depend on $\delta$\textup{)}. Suppose that $n > n_0(\omega,\delta)$ is sufficiently large, and let $f : \F_5^n \rightarrow [-1,1]$ be a function. Let $(\B^{(0)}_1,\B^{(0)}_2)$ be a quadratic factor of complexity $(d^{(0)}_1,d^{(0)}_2)$. Then there is $C = C(\delta,\omega,d^{(0)}_1,d^{(0)}_2)$ and a quadratic factor $(\B_1,\B_2)$ which refines $(\B^{(0)}_1,\B^{(0)}_2)$ and has complexity at most $(d,d)$, $d \leq C$, together with a decomposition
\[ f = f_1 + f_2 + f_3,\]
where
\[ f_1 := \E(f | \B_2),\]
\[ \Vert f_2 \Vert_2 \leq \delta\] and
\[ \Vert f_3 \Vert_{U^3} \leq 1/\omega(d).\]
\end{proposition}
\proof Apply the Koopman-von Neumann theorem iteratively, with parameters $\delta_i$, $i = 1,2,\dots$ to obtain quadratic factors $(\B_1^{(i)},\B_2^{(i)})$ with complexities at most $(C_i,C_i)$ such that
\begin{itemize}
\item $(\B_1^{(i)},\B_2^{(i)})$ is a refinement of $(\B_1^{(i-1)},\B_2^{(i-1)})$;
\item $\Vert f - \E(f | \B_2^{(i)}) \Vert_{U^3} \leq \delta_{i}$;
\item $C_{i}$ is bounded above in terms of $C_{i-1}$ and $\delta_{i}$.
\end{itemize}

Choose the sequence of $\delta_i$s such that $\delta_{i+1} \leq 1/\omega(C_i)$ for all $i$. Since $C_i$ is bounded above by a quantity depending only on $\delta_1,\dots,\delta_i$, this is certainly possible.

Now the energies $\Vert \E(f | \B_2^{(i)}) \Vert_2^2$ are non-decreasing, and are all bounded by 1. By the pigeonhole principle there is therefore some $i \leq \lceil \delta^{-2} \rceil$ such that 
\[ \Vert \E(f | \B_2^{(i+1)}) \Vert_2^2 - \Vert \E(f | \B_2^{(i)}) \Vert_2^2 \leq \delta^2.\]
For such an $i$, we may take for our decomposition
\[ f_1 := \E(f | \B_2^{(i)}),\]
\[ f_2 := \E(f | \B_2^{(i+1)}) - \E(f | \B_2^{(i)})\]
and 
\[ f_3 := f - \E(f | \B_2^{(i+1)}).\]
It follows from Pythagoras' Theorem that $\Vert f_2 \Vert_2 \leq \delta$, as required. \endproof

What is the point of the Koopman von Neumann and arithmetic regularity results, say for the $U^3$-norm? The answer is that they often reduce the study of general functions (say from the point of view of counting 4-term arithmetic progressions) to the study of projections $\E(f | \B)$ onto ``low-complexity'' quadratic factors. This, however, is of little consequence unless we can study those supposedly simple objects.

\begin{definition}[Rank of quadratic factors]
Suppose that $(\B_1,\B_2)$ is a quadratic factor of complexity $(d_1,d_2)$, being defined by $d_1$ linear forms $r_1^Tx,\dots, r_{d_1}^T x$ and $d_2$ pure quadratics $x^T M_1 x,\dots, x^T M_{d_2} x$. We say that $(\B_1,\B_2)$ has rank at least $r$ if
\[ \rk(\lambda_1 M_1 + \dots + \lambda_{d_2} M_{d_2}) \geq r\]
whenever $\lambda_1,\dots,\lambda_{d_2}$ are elements of $\F_5$, not all zero.
\end{definition}

When we are not concerned with bounds, it turns out that we may assume our quadratic factors have exceedingly large rank. We will see in the next lecture that factors with high rank are much easier to handle than factors with small rank.

\begin{lemma}[Making factors high-rank]\label{rank-red}
Let $\omega : \R_+ \rightarrow \R_+$ be an arbitrary growth function. Then there is another function $\tau = \tau_{\omega}$ with the following property. Let $(\B_1,\B_2)$ be a quadratic factor with complexity at most $(d_1,d_2)$. Then there is a refinement $(\B'_1,\B'_2)$ of $(\B_1,\B_2)$ with complexity at most $(d'_1,d_2)$, where $d'_1 \leq \tau(d_1,d_2)$, which has rank at least $\omega(d'_1 + d'_2)$.
\end{lemma}
\proof Suppose as usual that $(\B_1,\B_2)$ is described by $d_1$ linear functions $r_1^T x,\dots, r_{d_1}^Tx$ and $d_2$ ``pure quadratics'' $x^T M_1 x, \dots, x^T M_{d_2} x$. Suppose that $(\B_1,\B_2)$ does not have rank at least $\omega(d_1 + d_2)$.
Then there is some relation
\[ \rk(\lambda_1 M_1 + \dots + \lambda_{d_2} M_{d_2}) \leq \omega(d),\]
where we may assume without loss of generality that $\lambda_{d_2} = 1$. Let $s_1,\dots,s_k$, $k \leq \omega(d)$ be a basis for $\ker(U)^{\perp}$, where $U := \lambda_1 M_1 + \dots + \lambda_{d_2} M_{d_2}$, and let $(\B^{\dagger}_1,\B^{\dagger}_2)$ be the homogeneous quadratic factor defined by the linear forms $r_1^T x,\dots,r_{d_1}^T x, s_1^Tx, \dots, s_k^T x$ and the quadratic forms $x^T M_1 x,\dots, x^T M_{d_2-1} x$. It has complexity bounded by $(d^{\dagger}_1,d_2 - 1)$, where $d^{\dagger}_1 \leq d + \omega(d_1 + d_2)$.
The value of $x^T M_{d_2} x$ is determined by the values of the $x^T M_i x$, $i = 1,\dots,d_2-1$ together with the value of $x^T U x$. This in turn is determined by the coset of $\ker(U)$ that $x$ lies in, and hence by $s_1^Tx,\dots, s_k^T x$. It follows that $(\B^{\dagger}_1,\B^{\dagger}_2)$ refines $(\B_1,\B_2)$.

Now we ask whether $(\B^{\dagger}_1,\B^{\dagger}_2)$ has rank at most $\omega(d^{\dagger}_1 + d_2)$. If so, we refine again, obtaining a new factor $(\B^{\dagger \dagger}_1,\B^{\dagger \dagger}_2)$ with complexity bounded by $(d^{\dagger}_1 + \omega(d^{\dagger}_1 + d_2),d_2 - 2)$. This procedure can last no more than $d_2$ steps, however, since at each stage the number of pure quadratic phases is reduced by one. We may take $(\B'_1,\B'_2)$ to be the factor that we have when the procedure terminates.
\endproof

\begin{proposition}[Arithmetic regularity lemma for $U^3$, II]\label{gen-arith-reg-lem}
Let $\delta > 0$ be a parameter, and let $\omega_1,\omega_2 : \R_+ \rightarrow \R_+$ be arbitrary growth functions \textup{(}which may depend on $\delta$\textup{)}. Let $n > n_0(\delta,\omega_1,\omega_2)$ be sufficiently large, and let $f : \F_5^n \rightarrow [-1,1]$ be a function. Let $(\B^{(0)}_1,\B^{(0)}_2)$ be a quadratic factor of complexity $(d^{(0)}_1,d^{(0)}_2)$. Then there is a quadratic factor $(\B_1,\B_2)$ with the following properties:
\begin{enumerate}
\item $(\B_1,\B_2)$ refines $(\B^{(0)}_1,\B^{(0)}_2)$;
\item The complexity of $(\B_1,\B_2)$ is at most $(d_1,d_2)$, where \[ d_1,d_2 \leq C(\delta,\omega_1,\omega_2,d^{(0)}_1,d^{(0)}_2),\] for some fixed function $C$;
\item The rank of $(\B_1,\B_2)$ is at least $\omega_1(d_1 + d_2)$;
\item There is a decomposition $f = f_1 + f_2 + f_3$, where
\[ f_1 := \E(f | \B_2),\]
\[ \Vert f_2 \Vert_2 \leq \delta\] and
\[ \Vert f_3 \Vert_{U^3} \leq 1/\omega_2(d_1 + d_2).\]
\end{enumerate}
\end{proposition}
\emph{Remark.} The formulation is very similar to that in Proposition \ref{arith-reg-lem}, but we now insist that the factor $(\B_1,\B_2)$ be homogeneous, and also include a condition on its rank. The statement of Proposition \ref{gen-arith-reg-lem} will look complicated at first sight, but there is nothing much to be scared of. As always with complicated propositions, it is as well to attempt to formulate what has been proved in a somewhat looser, wordier way. Here is an attempt: 

\emph{Let $f$ be any function on $\F_5^n$. Then, up to an error which is small in $L^2$, we may write $f$ as a sum of a function which is measurable with respect to a bounded complexity quadratic factor, plus an error which is miniscule in $\Vert \cdot \Vert_{U^3}$. Furthermore we may insist that the rank of the quadratic factor is huge in comparison to its complexity.}

\proof 
Apply Proposition \ref{arith-reg-lem} to get a factor $(\B_1,\B_2)$ refining $(\B^{(0)}_1,\B^{(0)}_2)$, and a decomposition $f = f_1 + f_2 + f_3$ such that $f_1 = \E(f | \B_2)$, $\Vert f_2 \Vert_2 \leq \delta/2$ and $\Vert f_3 \Vert_{U^3} \leq 1/\omega_2(\tau(d_1, d_2) + d_2)$, where $(d_1,d_2)$ is an upper bound for the complexity of $(\B_1,\B_2)$ and $\tau = \tau_{\omega_1}$ is the function appearing in Lemma \ref{rank-red}. Using that lemma, we may refine $(\B_1,\B_2)$ to a quadratic factor $(\B'_1,\B'_2)$ with complexity at most $(d'_1,d'_2)$, where $d'_1 \leq \tau(d_1,d_2)$ and $d'_2 \leq d_2$, and with rank at least $\omega_1(d'_1 + d'_2)$. Define a new decomposition $f = f'_1 + f'_2 + f'_3$, where
\[ f'_1 := \E(f | \B'_2),\] 
\[ f'_2 := f_2 + \E(f | \B_2) - \E(f | \B'_2)\]
and $f'_3 = f_3$. Either this has the desired properties, or else we have 
\[ \Vert \E(f | \B_2) - \E(f | \B'_2) \Vert_2 \geq \delta/2.\] By Pythagoras' theorem this leads to the energy increment
\begin{equation}\label{en-inc-2} \Vert \E(f | \B'_2) \Vert_2^2 \geq \Vert \E(f | \B_2) \Vert_2^2 + \delta^2/4.\end{equation}
In this eventuality we apply Proposition \ref{arith-reg-lem} again, initialising with $(\B^{(0)}_1,\B^{(0)}_2) := (\B'_1,\B'_2)$. In view of the energy increment \eqref{en-inc-2}, we can only repeat this $\lceil 4/\delta^2 \rceil$ times before we reach a decomposition with the properties we desire.\endproof

\scriptsize
\emph{Further reading.} There is a wealth of directions to go in. Results of Koopman von Neumann type go back, implicitly, a long way. The \emph{name} was first given, by Tao and I, to a result in our paper \cite{green-tao-longprimeaps} on primes in AP. That result was somewhat different to the results here, but the method of proof (the energy increment strategy) is the same.

The arithmetic regularity lemma for the $U^3$-norm will be the subject of a forthcoming paper by Tao and I \cite{green-tao-tobewritten}. There is, of course, an analogous result for $U^2$-norm, and this was implicit in Bourgain \cite{bourgain-israel}. The proof there used the Fourier transform rather than the energy-increment strategy. A substantially more difficult (!) proof of the same result was given 15 years later by me \cite{green-gafa}; a number of applications were given there. The energy-increment proof of Proposition \ref{arith-reg-lem} seems at the moment to be the ``right'' way to think about these issues, and is essentially the approach taken in \cite{tao:ergodic}.

There are connections with regularity results for graphs and hypergraphs, the first result of this type being Szemer\'edi's regularity lemma \cite{sz-reg}. There are also parallels with results in ergodic theory such as \cite{bhk}. Perhaps it is best to refer the reader to the lectures by Kra and Tao at this school. The ICM article by Tao \cite{tao:icm2006} has many references and would represent a fine place to begin further investigations.
\normalsize

\section{Lecture 4}

Topics to be covered

\begin{itemize}
\item Working on a quadratic factor; the configuration space.
\item A theorem on progressions of length 4: an example of how to put all the ingredients together.
\end{itemize}

Our aim in this lecture is to prove the following theorem by using the machinery we have developed. Recall that we are writing $N := 5^n$.

\begin{theorem}[G.-Tao]\label{bhk-conjecture}
Let $\alpha,\epsilon > 0$ be real numbers. Then there is an $n_0 = n_0(\alpha,\epsilon)$ with the following property. Suppose that $n > n_0(\alpha,\epsilon)$, and that $A \subseteq \F_5^n$ is a set with density $\alpha$. Then there is some $d \neq 0$ such that $A$ contains at least $(\alpha^4 - \epsilon)N$ four-term arithmetic progressions with common difference $d$.
\end{theorem}
\emph{Remarks.} It is easy to see that one cannot replace $\alpha^4$ by anything larger, by considering a random set of density $\alpha$. This theorem has, as a consequence, a version of Szemer\'edi's theorem for progressions of length four in finite fields, namely $r_4(\F_5^n) = o(N)$. The theorem is a finite field version of a conjecture of Bergelson, Host and Kra. Rather bizarrely at first sight, this result does not generalise to progressions longer than four.

Now in the last lecture we worked rather hard in order to show that, in various senses, the study of an arbitrary function $f : \F_5^n \rightarrow [-1,1]$ can be reduced to the study of a $\B_2$-measurable function $\E(f | \B_2)$, where $(\B_1,\B_2)$ is a  quadratic factor with ``bounded complexity'' and high rank. To make use of this, we need to be able to understand $\B_2$-measurable functions. At the very least, we are going to want to know about the size of the atoms in $\B_2$ and, for any four atoms, the number of four-term progressions spanned by those atoms. It turns out that the ``high-rank'' assumption allows us to simply compute these quantities using Fourier analysis.

Suppose, throughout this lecture, that $(\B_1,\B_2)$ is a  quadratic factor defined by $d_1$ linear forms $r_j^Tx$ and $d_2$ pure quadratics $x^T M_j x$. (Recall that $\B_1$ is the $\sigma$-algebra generated by the linear functions, and $\B_2$ is the $\sigma$-algebra generated by the linear \emph{and} quadratic functions.) We will always suppose (as we clearly may) that the vectors $r_j$ are linearly independent. 

To understand $\B_2$-measurable functions, that is to say functions which are constant on atoms of $\B_2$ (or alternatively functions which have the form $\E(f|\B_2)$), it is helpful to work in \emph{configuration space} $\F_5^{d_1} \times \F_5^{d_2}$. We write $\Gamma: \F_5^n \to \F_5^{d_1}$ and $\Phi: \F_5^n \to \F_5^{d_2}$ for the maps
$\Gamma(x) := (r_1^Tx,\dots,r_{d_1}^Tx)$ and $\Phi(x) := (x^T M_1 x,\dots,x^T M_{d_2}x)$.

\begin{lemma}[Size of atoms]\label{at-lem} Suppose that $(\B_1,\B_2)$ has rank at least $r$.
Let $(a,b) \in \F_5^{d_1} \times \F_5^{d_2}$. Then the probability that a randomly chosen $x \in \F_5^n$ has $\Gamma(x) = a$ and $\Phi(x) = b$ is $5^{-d_1 - d_2} + O(5^{- r/2})$.
\end{lemma}
\emph{Remark.} In this lemma and the next, the probabilistic language is present only to avoid normalising factors of $N = 5^n$. This is really a statement about the \emph{number} of $x$ with $\Gamma(x) = a$, $\Phi(x) = b$.

\proof The quantity in question is given by
\[ 5^{-d_1-d_2}\E_x \prod_{i=1}^{d_1} \big(\sum_{\mu_i \in \F_5} \omega^{\mu_i(r_i^T x - a_j)}\big)\prod_{j=1}^{d_2} \big(\sum_{\lambda_j \in \F_5} \omega^{\lambda_j(x^T M_j x - b_j)}\big),\] which rearranges as
\begin{equation}\label{star-2} 5^{-d_1-d_2}\sum_{\mu_i,\lambda_j} \omega^{-\lambda_1 b_1 - \dots - \lambda_{d_2}b_{d_2} - \mu_1 a_1 - \dots - \mu_{d_1} a_{d_1}}\E_x \omega^{x^T(\lambda_1 M_1 + \dots + \lambda_{d_2}M_{d_2})x + (\mu_1 r_1 + \dots + \mu_{d_1}r_{d_1})^T x} .\end{equation}
Now the rank of $(\B_1,\B_2)$ is at least $r$, which means that
\[ \rk(\lambda_1 M_1 + \dots + \lambda_{d_2}M_{d_2}) \geq r.\]
In view of the Gauss sum estimate, Lemma \ref{gauss-sums}, this means that every term in \eqref{star-2} in which the $\lambda_i$ are not all zero is bounded by $5^{-d_1 -d_2 -r/2}$. Of the terms with $\lambda_1 = \dots = \lambda_{d_2} = 0$, the linear independence of the $r_i$ guarantees that the only term which does not vanish is that with $\mu_1 = \dots = \mu_{d_1} = 0$. The result follows immediately.\endproof

\begin{lemma}[4-term progressions]\label{4ap-lem} Suppose that $(\B_1,\B_2)$ has rank at least $r$.
Suppose that $(a^{(1)},b^{(1)}),\dots,(a^{(4)},b^{(4)}) \in \F_5^{d_1} \times \F_5^{d_2}$. Suppose that a 4-term progression $(x,x+d,x+2d,x+3d) \in (\F_5^n)^4$ is chosen at random. If \begin{equation}\label{constraint-1} \mbox{$a^{(1)},a^{(2)},a^{(3)},a^{(4)}$ are in arithmetic progression}\end{equation} and \begin{equation}\label{constraint-2} b^{(1)} - 3b^{(2)} + 3b^{(3)} + b^{(4)} = 0\end{equation} then the probability that $\Gamma(x + id) = a^{(i)}$, $\Phi(x + id) = b^{(i)}$ for $i = 1,2,3,4$ is $5^{-2d_1 - 3d_2} + O(5^{-r/2})$. Otherwise, it is zero.
\end{lemma}
\proof The important thing to appreciate here is that four elements in different atoms of $\B_2$ can only lie in arithmetic progression if the two constraints \eqref{constraint-1} and \eqref{constraint-2} are satisfied. Furthermore these are the only relevant constraints, in that if they are satisfied (and if the factor $(\B_1,\B_2)$ has large rank) then we can accurately count the number of four-term progressions involving those atoms.

The necessity of the constraints \eqref{constraint-1} and \eqref{constraint-2} is easy. If $(x,x+d,x+2d,x+3d)$ is an arithmetic progression, we need only observe that $\Gamma(x),\Gamma(x+d),\Gamma(x+2d),\Gamma(x+3d)$ are also in arithmetic progression, and that $\Phi(x) - 3\Phi(x+d) + 3\Phi(x+2d) - \Phi(x+3d) = 0$.

To obtain the statement about probability, we proceed in the same manner as in Lemma \ref{at-lem}. The notation here is, however, somewhat fearsome. We start with the observation that the probability in question is
\[ 5^{-4d_1 - 4d_2}\E_{x,d} \prod_{l=1}^4\prod_{i=1}^{d_1} \big( \sum_{\mu_i^{(l)} \in \F_5} \omega^{\mu^{(l)}_i (r_i^T (x + ld) - a_i^{(l)})} \big) \prod_{j=1}^{d_2} \big( \sum_{\lambda_j^{(l)} \in \F_5} \omega^{\lambda_j^{(l)} ((x + ld)^T M_j (x + ld) - b_j^{(l)})} \big),\] and then swap the order of summation to rearrange as
\begin{equation}\label{star-3} 5^{-4d_1-4d_2}\sum_{\mu^{(l)}_i,\lambda^{(l)}_j \in \F_5} \E_{x,d}\omega^{x^T Px + 2x^T Q d + d^T R d + u^T x + v^T d - w},\end{equation} where 
\[ P = P(\lambda) = \sum_{j = 1}^{d_2}(\lambda^{(1)}_j + \lambda^{(2)}_j + \lambda^{(3)}_j + \lambda^{(4)}_j) M_j,\]
\[ Q = Q(\lambda) = \sum_{j = 1}^{d_2}(\lambda^{(1)}_j + 2\lambda^{(2)}_j + 3\lambda^{(3)}_j + 4\lambda^{(4)}_j) M_j,\]
\[ R = R(\lambda) = \sum_{j = 1}^{d_2}(\lambda^{(1)}_j + 4\lambda^{(2)}_j + 9\lambda^{(3)}_j + 16\lambda^{(4)}_j) M_j,\]
\[ u = u(\mu) = \sum_{i = 1}^{d_1} (\mu^{(1)}_i + \mu^{(2)}_i + \mu^{(3)}_i + \mu^{(4)}_i) r_i,\]
\[ v = v(\mu) = \sum_{i=1}^{d_1}(\mu^{(1)}_i + 2\mu^{(2)}_i + 3\mu^{(3)}_i + 4\mu^{(4)}_i) r_i\]
and
\[ w = w(\mu,\lambda) = \sum_{l=1}^4 \sum_{i=1}^{d_1} \mu^{(l)}_i a^{(l)}_i + \sum_{l = 1}^4 \sum_{j = 1}^{d_2} \lambda^{(l)}_j b^{(l)}_j.\]
We use Lemma \ref{gauss-sums} repeatedly. By fixing either $x$ or $d$, we see that the inner sum in \eqref{star-3} (that is, the expectation over $x,d$) is $O(5^{-r/2})$ unless \begin{equation}\label{cond-1} \lambda^{(1)}_j + \lambda^{(2)}_j + \lambda^{(3)}_j + \lambda^{(4)}_j = \lambda^{(1)}_j + 4\lambda^{(2)}_j + 9\lambda^{(3)}_j + 16\lambda^{(4)}_j = 0,\end{equation} in which case certainly $P = R = 0$. In this case, the inner sum is a rather purer-looking
\begin{equation}\label{new-inner-sum} \E_{x,d} \omega^{x^T Q d + u^T x + v^T d - w}.\end{equation}
For fixed $d$, this is zero unless $Qd + u = 0$. If $\lambda^{(1)}_j + 2\lambda^{(2)}_j + 3\lambda^{(3)}_j + 4\lambda^{(4)}_j \neq 0$ then, since $\rk(Q) \geq r$, this cannot happen for more than $5^{-r}$ of all $d$, and \eqref{new-inner-sum} is bounded by $5^{-r}$. If on the other hand \begin{equation}\label{cond-2} \lambda^{(1)}_j + 2\lambda^{(2)}_j + 3\lambda^{(3)}_j + 4\lambda^{(4)}_j = 0\end{equation} then \eqref{new-inner-sum} further reduces to
\[ \E_{x,d} \omega^{u^Tx + v^T d - w},\] which clearly vanishes unless 
\begin{equation}\label{cond-3} \mu^{(1)}_i + \mu^{(2)}_i + \mu^{(3)}_i + \mu^{(4)}_i = \mu^{(1)}_i + 2\mu^{(2)}_i + 3\mu^{(3)}_i + 4\mu^{(4)}_i = 0.\end{equation}
We have shown that the inner sum in \eqref{star-3} is $O(5^{-r/2})$ unless the five linear conditions \eqref{cond-1},\eqref{cond-2},\eqref{cond-3} are satisfied. The total contribution to \eqref{star-3} from cases where one of these five conditions is not satisfied is therefore $O(5^{-r/2})$. The total contribution from cases when the five conditions \emph{are} satisfied is
\[ 5^{-4d_1 - 4d_2} \sum_{l = 1}^4 \sum_{\mu^{(l)}_i,\lambda^{(l)}_j} \omega^{-w(\mu,\lambda)}.\]
Since the $a^{(i)}$ are in arithmetic progression and the $b^{(i)}$ satisfy $b^{(1)} - 3b^{(2)} + 3b^{(3)} - b^{(4)}$, it is easy to check that $w(\mu,\lambda) = 0$ when the five conditions are satisfied. It remains only to note that, of the $5^{4d_1 + 4d_2}$ choices for $\mu,\lambda$, the five conditions are satisfied for $5^{2d_1 + d_2}$ of them.\endproof

If $f : \F_5^n \rightarrow \mathbb{C}$ is a $\B$-measurable function then we write $\f : \F_5^{d_1} \times \F_5^{d_2} \rightarrow \mathbb{C}$ for the function which satisfies
\[ f(x) = \f(\Gamma(x), \Phi(x))\] for all $x \in \F_5^n$. We will adopt this convention of using bold letters to denote functions on configuration space for the rest of these lectures without further comment.

We are now in a position to prove Theorem \ref{bhk-conjecture}.

\emph{Proof of Theorem \ref{bhk-conjecture}.} Recall that $A \subseteq \F_5^n$ is a set with density $\alpha$. Apply Proposition \ref{gen-arith-reg-lem} to find a quadratic factor $(\B_1,\B_2)$ with complexity $(d_1,d_2)$, $d_i \leq d_0(\alpha,\epsilon)$ and rank $r$ satisfying (say)
\[ r \geq 100(\log(1/\epsilon) + \log(1/\alpha) + d_1 + d_2)\]together with a decomposition $1_A = f_1 + f_2 + f_3$ such that $f_1 = \E(1_A | \B_2)$, $\Vert f_2 \Vert_2 \leq \delta$ and $\Vert f_3 \Vert_{U^3} \leq 1/\omega(d_1 + d_2)$. The parameter $\delta$ and the growth function $\omega$ will be specified as the proof unfolds, but will depend only on $\alpha$ and $\epsilon$. 

Let $r_1^T x,\dots,r^T_{d_1}x$ be the linear functions involved in $\B_1$, and let $H := \langle r_1,\dots,r_{d_1}\rangle^T$. Let $1_H$ be the characteristic function of $H$, and let $\mu_H$ be the normalised measure on $H$, thus $\mu_H := 1_H/\E1_H$. We are going to prove that 
\begin{equation}\label{to-prove} \E_{x,d} 1_A(x)1_A(x+d)1_A(x+2d)1_A(x+3d) \mu_H(d) \geq \alpha^4 - \epsilon,\end{equation} which clearly implies the theorem (for some $d \in H$). To do this, we split the left-hand-side of \eqref{to-prove} into 81 parts by substituting $1_A = f_1 + f_2 + f_3$.

\emph{Claim 1.} The contribution from any of the 65 terms which contain $f_2$ is no more $\epsilon/200$. 

\proof Suppose that the term is
\begin{equation}\label{to-est-1} \E_{x,d} g_1(x) g_2(x+d)g_3(x+2d) g_4(x+3d)\mu_H(d),\end{equation} where $g_1 = f_2$ (the proofs of the other cases are very similar). Set $F(x) := \E_d g_2(x+d)g_3(x+2d)g_4(x+3d) \mu_H(d)$, and observe that $\Vert F \Vert_{\infty} \leq 1$. It follows that 
\[ |\E_{x,d} g_1(x) g_2(x+d)g_3(x+2d) g_4(x+3d)\mu_H(d)| \leq |\E_{x} g_1(x) F(x)| \leq \Vert f_2 \Vert_1 \leq \Vert f_2 \Vert_2.\] This proves the claim provided that $\delta \leq \epsilon/200$.

\emph{Claim 2.} The contribution from any of the 65 terms which contain $f_3$ is no more than $\epsilon/200$.

\proof Suppose that the term is
\begin{equation}\label{to-est-2} \E_{x,d} g_1(x) g_2(x+d)g_3(x+2d) g_4(x+3d)\mu_H(d),\end{equation} where $g_1 = f_3$ (the proofs of the other cases are very similar). We have
\[ 1_H(d) = \sum_t 1_{t+H}(x+2d) 1_{t+H}(x +d), \]
where the sum is over all cosets $t + H$ of $H$ in $\F_5^n$. By the generalised von Neumann theorem (Proposition \ref{4ap-gvn}), we have 
\[ |\E_{x,d} g_1(x) g_2(x+d) 1_{t+H}(x+d) g_3(x+2d) 1_{t+H}(x+2d) g_4(x+3d)| \leq \Vert f_3 \Vert_{U^3} \leq 1/\omega(d_1 +d_2)\] for each $t$. It follows that \eqref{to-est-2} is no more than $5^{2d_1}/\omega(d_1 + d_2)$, which proves the claim provided that $\omega(t) \geq 5^{t+4}/\epsilon$.

\emph{Remarks.} Note carefully that for Claim 2 to follow we required the regularity parameter $\omega(t)$ to be \emph{exponential} in $t$, rather than (say) polynomial. This is why the full arithmetic regularity lemma is required, rather than just the Koopman-von Neumann theorem.

These two claims account for 80 of the 81 terms into which we have decomposed the left-hand side of \eqref{to-prove}. To finish the argument, it suffices to show that
\begin{equation}\label{to-prove-7} \E_{x,d} f_1(x)f_1(x+d)f_1(x+2d)f_1(x+3d) \mu_H(d) \geq \alpha^4 - \epsilon/2.\end{equation} 
Now $f_1$ is (by definition) constant on atoms of $\B_2$. Recall that these atoms are indexed by the configuration space $\F_5^{d_1} \times \F_5^{d_2}$, and that we write $\f_1(a,b)$ for the value of $f_1$ on the atom indexed by $(a,b)$. 

\emph{Claim 3.} We have
\begin{equation}\label{dens}
\E_{(a,b) \in \F_5^{d_1} \times \F_5^{d_2}} \f_1(a,b) = \alpha(1 + O(5^{2d_1 + 2d_2 - r/2})).
\end{equation}
\proof Note that the result would be trivial (and would hold without the $O$-term) if all the atoms of $\B_2$ had \emph{exactly} the same size. Now recall that Lemma \ref{at-lem} gives an approximate version of this statement. We leave the slightly tedious details to the reader.

\emph{Claim 4.} We have 
\begin{align*}
\E_{x,d} & f_1(x)f_1(x+d)f_1(x+2d)f_1(x+3d) \mu_H(d) \\ &= \E_{\substack{a \in \F_5^{d_1} , b^{(1)},\dots,b^{(4)} \in \F_5^{d_2} \\ b^{(1)} - 3b^{(2)} + 3b^{(3)} - b^{(4)} = 0}} \f_1(a,b^{(1)})\f_1(a,b^{(2)})\f_1(a,b^{(3)})\f_1(a,b^{(4)}) + O(5^{2d_2 +3d_2 - r/2}).
\end{align*}
\proof Condition on the quadruple $(a^{(1)},b^{(1)}),\dots,(a^{(4)},b^{(4)})$ of atoms containing $(x,x+d,x+2d,x+3d)$. The constraint that $d  \in H$ is equivalent to $a^{(1)} = a^{(2)} = a^{(3)} = a^{(4)} = a$, say. By Lemma \ref{4ap-lem}, we must also have $b^{(1)} - 3b^{(2)} + 3b^{(3)} - b^{(4)} = 0$. Invoking that same lemma, we have
\begin{align*}
& \E_{x,d} f_1(x)f_1(x+d)f_1(x+2d)f_1(x+3d) 1_H(d) \\
&= (5^{-2d_1 - 3d_2} + O(5^{-r/2}))\sum_{a \in \F_5^{d_1}} \sum_{\substack{b^{(1)},\dots,b^{(4)} \in \F_5^{d_2} \\ b^{(1)} - 3b^{(2)} + 3b^{(3)} - b^{(4)} = 0}} \f_1(a,b^{(1)})\f_1(a,b^{(2)})\f_1(a,b^{(3)})\f_1(a,b^{(4)}) 
\end{align*}
Normalising, we obtain the stated result.

Now the rank $r$ was chosen very large ($r > 100(\log(1/\epsilon) + \log(1/\alpha) + d_1 + d_2)$). All we need do to establish \eqref{to-prove-7}, then, is prove the inequality

\begin{equation}\label{to-prove-8} \E_{\substack{a \in \F_5^{d_1}, b^{(1)},\dots,b^{(4)} \in \F_5^{d_2} \\ b^{(1)} - 3b^{(2)} + 3b^{(3)} - b^{(4)} = 0}} \f_1(a,b^{(1)})\f_1(a,b^{(2)})\f_1(a,b^{(3)})\f_1(a,b^{(4)}) \geq \big(\E_{(a,b) \in \F_5^{d_1} \times \F_5^{d_2}} \f_1(a,b)\big)^4.\end{equation}

Noting that the left-hand side is 
\[ \E_{a \in \F_5^{d_1}} \E_{x \in \F_5^{d_2}} \big( \E_{\substack{b,b' \in \F_5^{d_2}\\ b - 3b' = x}} \f_1(a,b) \f_1(a,b') \big)^2,\]
this follows from two applications of the Cauchy-Schwarz inequality.

Alternatively, it is amusing to give an interpretation in terms of the Fourier transform. The left-hand side of \eqref{to-prove-8} is

\begin{equation}\label{four-exp} \E_{a \in \F_5^{d_1}} \sum_{r \in \widehat{\F_5^{d_2}}} |\widetilde{\f_1}(a,r)|^2 |\widetilde{\f_1}(a,-3r)|^2.\end{equation}
In this expression the tilde denotes Fourier transform in the second variable, which was called $b$ in \eqref{to-prove-8}.

A lower bound for \eqref{four-exp} comes from ignoring all terms except those with $r = 0$, yielding
\[ \E_{a \in \F_5^n} |\widetilde{\f_1}(a,0)|^4 = \E_{a \in \F_5^n} |\E_{b \in \F_5^{d_2}} \f_1(a,b)|^4.\] 

The result now follows from H\"older's inequality.

A more interesting application of these partial Fourier transforms may be found in \cite{green-tao-ffszem}.

\section{Lecture 5}
Topics to be covered

\begin{itemize}
\item An introduction to the theory on $\Z/N\Z$.
\end{itemize}

For simplicity I will assume that $N$ is a large prime.

I am only scheduled to give four lectures at the school. These notes are here for two reasons: firstly, it is possible that I will finish the material from the first four lectures early. More importantly, it is the theory on the group $\Z/N\Z$ that is of most interest for applications in number theory, and it would be remiss of me to not at least point the reader in directions where she may learn more.

Note that the theory on $\Z/N\Z$ is actually rather richer than for an arbitrary abelian group $G$, because we have been able to pursue analogies with ergodic theory. This is concerned with $\Z$-actions, and $\Z/N\Z$ is the finite abelian group which most closely models $\Z$.

One way of motivating the theory is to try and take what we know for $\F_5^n$ and attempt to adapt it to $\Z/N\Z$. Let us note that the basic definitions of Gowers norms and the basic generalised von Neumann theorems of Lecture 1 go over essentially unchanged to $\Z/N\Z$. The first stumbling block comes at the point where we ask for a conjectural analogue of Proposition \ref{inv-thm-1}. A first guess might be:

\begin{conjecture}\label{wrong-conj} Suppose that $f : \Z/N\Z \rightarrow [-1,1]$ is a function with $\Vert f \Vert_{U^3} \geq \delta$. Then there are $r,s \in \Z/N\Z$ such that
\[ |\E_{x \in \Z/N\Z} f(x) e\big( \frac{rx^2 + sx}{N} \big)| \gg_{\delta} 1.\]
\end{conjecture}
\emph{Remark.} As usual in analytic number theory we have written $e(\theta) := e^{2\pi i \theta}$.

It turns out that this conjecture is false. One example of a function on $\Z/N\Z$ which has large $U^3$-norm, but does not correlate with a quadratic form $e(rx^2 + sx/N)$, is a quadratic $e(\theta x^2)$ where $\theta \not\approx r/N$. Such a quadratic is most naturally defined on $\Z$, but by restricting its domain to $\{1,\dots,N\}$ one obtains a function which can be defined on $\Z/N\Z$. Another example is a ``bracket quadratic'' such as $e(\theta_1 x \{\theta_2 x\})$, where $\{t\}$ denotes the fractional part of $t$. The second of these counterexamples is somehow more serious, but it is also rather harder to see that this rather exotic function \emph{does} provide a counterexample to Conjecture \ref{wrong-conj}. For a brief discussion see \cite[\S 6]{green-icm}, and for more detail see \cite{green-tao-u3inverse}.

If the only obvious generalisation of Proposition \ref{inv-thm-1} is wrong, how should we proceed? It turns out that a hint is given to us by the quantitatively stronger form of the inverse theorem for the $U^3$-norm on $\F_5^n$, namely Proposition \ref{inv-thm-2}. We are not concerned with quantitative issues here, so let us state a weak consequence of that result. This is actually a trivial consequence of Proposition \ref{inv-thm-1}, too.

\begin{proposition}[Inverse result for $U^3$-norm on $\F_5^n$, III]\label{inv-thm-3} Suppose that $f : \F_5^n \rightarrow [-1,1]$ is a function with $\Vert f \Vert_{U^3} \geq \delta$. Then there is a subspace $H \leq \F_5^n$ with $\codim H \ll_{\delta} 1$, a matrix $M \in \GL_n(\F_5)$ and a vector $r \in \F_5^n$ such that
\[ |\E_x f(x)1_H(x) \omega^{x^T M x + r^T x}| \gg_{\delta} 1.\]
\end{proposition}
\emph{Remark.} It is not too hard to show that this is \emph{equivalent} to Proposition \ref{inv-thm-1}: we leave this as an exercise to the reader.

Let us try and generalise this result. There are two objects which do not obviously transfer to $\Z/N\Z$: the notion of \emph{subspace}, and (implicitly) the notion of \emph{quadratic form}. It turns out that the second notion can be sensibly formulated for functions defined on any set.

\begin{definition}[Quadratic forms]
Let $S$ be a set in some abelian group, and let $\psi : S \rightarrow \R/\Z$ be a function. We say that $\psi$ is a quadratic form if the second derivative
\[ \psi''(h_1,h_2) := \psi(x+h_1+h_2) - \psi(x + h_1) - \psi(x + h_2) + \psi(x)\]
is well-defined, that is to say if this definition does not depend on $x$ whenever $x,x+h_1,x+h_2,x+h_1+h_2 \in S$.
\end{definition}

Whilst the notion of subspace is rather vacuous in $\Z/N\Z$, there is a plentiful supply of \emph{approximate} subspaces. These are more usually called Bohr sets.

\begin{definition}[Approximate subspaces/Bohr sets] Let $R =\{r_1,\dots,r_k\} \subseteq \Z/N\Z$ and let $\epsilon > 0$. Then we write
\[ B(R,\epsilon) := \{ x \in \Z/N\Z : |e(rx/N) - 1| \leq \epsilon\}.\]
This is called the Bohr set with width $\epsilon$ corresponding to frequency set $R$.
\end{definition}
The set $R$ should actually be thought of as a set of characters on $\Z/N\Z$, each value $r$ corresponding to the character $x \mapsto e(rx/N)$. Once thought of in this way, it is easy to see how Bohr sets can be defined on any finite abelian group $G$. Bohr sets on $\F_5^n$ do not depend very seriously on the width parameter $\epsilon$, and certainly for $\epsilon < 1/10$ (say) they are just vector subspaces.

There is a lot to say about Bohr sets, and much information may be found in \cite{tao-vu-book}. See also \cite{green-fin-field}, where there is a discussion of the place of Bohr sets in the transition from finite field models to $\Z/N\Z$ in various settings. We caution the reader that there are certain technicalities associated with the study of Bohr sets in additive combinatorics, most particularly the need to consider \emph{regular} Bohr sets (ones that ``behave well at the edges''). In this brief overview we will say nothing more about these technicalities, other than that most of them were overcome in a seminal paper of Bourgain \cite{Bou}.

To return to the point, we may now state Theorem 2.7 (i) of \cite{green-tao-u3inverse}, which is an inverse theorem for the $U^3$-norm on $\Z/N\Z$. In the light of the above discussion, the reader will see that it is a natural generalisation of Proposition \ref{inv-thm-3}.

\begin{proposition}[Inverse theorem for the $U^3$-norm on $\Z/N\Z$, I]\label{zn-inv}
Suppose that $f : \Z/N\Z \rightarrow [-1,1]$ is a function and that $\Vert f \Vert_{U^3} \geq \delta$. Then there is a set $R \subseteq \Z/N\Z$, $|R| \ll_{\delta} 1$, a parameter $\epsilon \gg_{\delta} 1$ such that the Bohr set $B := B(R,\epsilon)$ is regular, some $y \in \Z/N\Z$ and a quadratic form $\psi : y + B \rightarrow \R/\Z$ such that
\begin{equation}\label{eq497} |\E_x f(x) 1_{y + B}(x) e(\psi(x))| \gg_{\delta} 1.\end{equation}
\end{proposition}

It turns out that result is necessary and sufficient, that is to say if \eqref{eq497} is satisfied then $\Vert f \Vert_{U^3}$ is large. See \cite{green-tao-u3inverse}, Thm 2.7 (ii) (note that this is the only point at which the regularity of $B(R,\epsilon)$ is relevant). This is, at first sight, a very unsatisfactory state of affairs: we have a theorem which gives a necessary and sufficient condition for a natural problem which interests us, yet the theorem is somewhat inelegant and difficult to state. 

Our subject being in some sense an extension of the work of Hardy and Littlewood, one should perhaps recall at this point Hardy's view that there is ``no permanent place in the world for ugly mathematics''. 

With this in mind we observe that although Proposition \ref{zn-inv} is necessary and sufficient, it need not be the \emph{only} necessary and sufficient condition. In what follows we will be rather vague. Write  $\mathcal{Q} = \mathcal{Q}(\delta)$ for the collection of all ``quadratic obstructions'' of the form $1_{y + B}(x) e(\psi(x))$, where $B,\psi$ are as above. Any other collection $\mathcal{Q}'$ with the property that anything in $\mathcal{Q}$ is approximately a linear combination of elements in $\mathcal{Q}'$, and vice versa, will also be a necessary and sufficient collection of quadratic obstructions for $\Z/N\Z$.

It turns out that there is a very natural choice for $\mathcal{Q}'$, the collection of $2$-\emph{step nilsequences}. The idea that we should look at these objects came to us from ergodic theory -- there will be much more on this in the lectures of Bryna Kra at the school.

Let $G$ be a connected, simply-connected $2$-step nilpotent Lie group over $\R$ and let $\Gamma \leq G$ be a discrete, cocompact submanifold. The quotient $G/\Gamma$ is called a $2$-step nilmanifold. For the sake of illustration, we recommend that the reader take

\[ G := \begin{pmatrix} 1 & \R & \R \\ 0 & 1 & \R \\ 0 & 0 & 1\end{pmatrix}, \qquad \Gamma := \begin{pmatrix} 1 & \Z & \Z \\ 0 & 1 & \Z \\ 0 & 0 & 1\end{pmatrix},\]
in which case $G/\Gamma$ is a 3-dimensional compact manifold called the \emph{Heisenberg nilmanifold}.

Let $g \in G$ and $x \in G/\Gamma$ be arbitrary. The element $g$ induces a continuous map $T_g : G/\Gamma \rightarrow G/\Gamma$ by multiplication on the left. Any sequence of the form $(F(T_g^n \cdot x))_{n \in \mathbb{N}}$, where $F : G/\Gamma \rightarrow [-1,1]$ is continuous, is called a 2-step nilsequence. It turns out that the collection of $2$-step nilsequences can play the r\^ole of $\mathcal{Q}'$ as discussed above. The following is proved in \cite{green-tao-u3inverse}, Thm. 12.8.

\begin{proposition}[Inverse Theorem for the $U^3$-norm on $\Z/N\Z$, II] Let $f : \Z/N\Z \rightarrow [-1,1]$ be a function, and suppose that $\Vert f \Vert_{U^3} \geq \delta$. Then there is a $2$-step nilsequence $(F(T_g^n \cdot x))_{n \in \mathbb{N}}$ with complexity $\ll_{\delta} 1$ such that 
\[ |\E_{n \leq N} f(n) F(T_g^n \cdot x)| \gg_{\delta} 1.\]
If, conversely, $f$ correlates with a $2$-step nilsequence of bounded complexity then the $\Vert \cdot \Vert_{U^3}$-norm of $f$ is large.
\end{proposition}
We have not defined the \emph{complexity} of a nilsequence. It is some number associated to $(F(T_g^n \cdot x))_{n \in \mathbb{N}}$, which bounds both the dimension of the underlying nilmanifold $G/\Gamma$, and also the Lipschitz constant of $F$ with respect to some sensible metric. There is no canonical way of defining the complexity, but this is not important for the theory.

We do not attempt to explain why this collection $\mathcal{Q}'$ of $2$-step nilsequences is ``equivalent'' to the collection $\mathcal{Q}$ used in Proposition \ref{zn-inv}. Detailed technical discussions may be found in \cite{green-tao-u3inverse,green-tao-u3mobius}. A short calculation involving the Heisenberg example, showing how a $2$-step nilsequence on it resembles a quadratic form on a Bohr set, is given in \cite{green-icm}.


\begin{thebibliography}{99}

\bibitem{balog}
A.~Balog and E.~Szemer\'edi, \emph{A statistical theorem of set addition},
Combinatorica, \textbf{14} (1994), 263--268.

\bibitem{bhk}
V.~Bergelson, B.~Host and B.~Kra, \emph{Multiple recurrence and nilsequences}, with an appendix by I.Z. Ruzsa, Invent. Math. \textbf{160} (2005), no. 2, 261--303.

\bibitem{bourgain-israel} J.~Bourgain, \emph{A Szemer\'edi-type theorem for sets of positive density in $\R^k$,} Israel J. Math. \textbf{54} (1986), no. 3, 307--316.

\bibitem{Bou} \bysame, \emph{On triples in arithmetic progression}, GAFA \textbf{9} (1999), no. 5, 968--984.

\bibitem{chang-er}
M.~C.~Chang, \emph{On problems of Erd\H{o}s and Rudin}, J. Funct. Anal. \textbf{207} (2004), 444--460.

\bibitem{croot-ap-lims} E.~Croot, \emph{The minimal number of 3-term arithmetic progressions modulo a prime converges to a limit,} to appear in Canadian Math. Bull.

\bibitem{gowers-example} W.~T.~Gowers, \emph{Lower bounds of tower type for Szemer\'edi's uniformity lemma,} GAFA \textbf{7} (1997), no. 2, 322-337.

\bibitem{gowers-4aps} \bysame, \emph{A new proof of Szemer\'edi's theorem for progressions of length four}, GAFA \textbf{8} (1998), no. 3, 529--551.


\bibitem{gowers-longaps} \bysame, \emph{A new proof of Szemer\'edi's theorem}, GAFA \textbf{11} (2001), 465--588.

\bibitem{green-roth-primes} B.~J.~Green, \emph{Roth's theorem in the primes,} Annals of Math. \textbf{161} (2005), no. 3, 1609--1636.

\bibitem{green-gafa} \bysame, \emph{A Szemer\'edi-type regularity lemma in abelian groups,} GAFA. \textbf{15} (2005), no. 2, 340--376.

\bibitem{green-fin-field} \bysame, \emph{Finite field models in additive combinatorics} in Surveys in Combinatorics 2005, London Math. Soc. Lecture Notes \textbf{327}, 1--27.

\bibitem{green-mit} \bysame, \emph{Edinburgh-MIT lecture notes on Freiman's theorem}, available at\\
\texttt{http://www.dpmms.cam.ac.uk/$\widetilde{\;}$bjg23}

\bibitem{green-bkt} \bysame, \emph{Notes on the Bourgain-Katz-Tao theorem,} 
available at\\
\texttt{http://www.dpmms.cam.ac.uk/$\widetilde{\;}$bjg23}

\bibitem{green-icm} \bysame, \emph{Generalising the Hardy-Littlewood method for primes,} Proceddings of the International Congress of Mathematicians, Madrid 2006, Vol. 2.

\bibitem{green-konyagin} B.~J.~Green and S.~Konyagin, \emph{On the Littlewood problem modulo a prime,} to appear in Canadian J. Math.

\bibitem{green-tao-longprimeaps} B.~J.~Green and T.~C.~Tao, \emph{The primes contain arbitrarily long arithmetic progressions}, to appear in Annals of Math.

\bibitem{green-tao-u3inverse} \bysame, \emph{An inverse theorem for the Gowers $U^3$-norm, with applications,} to appear in Proc. Edinburgh Math. Soc.

\bibitem{green-tao-ffszem} \bysame, \emph{A new bound for Szemer\'edi's theorem in finite field geometries, for progressions of length 4}, preprint.

\bibitem{green-tao-u3mobius} \bysame, \emph{Quadratic uniformity of the M\"obius function,} preprint.

\bibitem{green-tao-linearprimes} \bysame, \emph{Linear equations in primes,} preprint.

\bibitem{green-tao-tobewritten} \bysame, \emph{Arithmetic regularity lemmas,} to be written.

\bibitem{heath-brown-3aps} D.~R.~Heath-Brown, \emph{Integer sets containing no arithmetic progressions,} J. London Math. Soc. \textbf{35} (1987), 385--394.

\bibitem{host-kra}
B.~Host and B.~Kra, \emph{Non-conventional ergodic averages and nilmanifolds,} Ann. Math. \textbf{161} (2005), no. 1, 397--488.

\bibitem{plun}
H.~Pl\"unnecke, \emph{Eigenschaften un Absch\"atzungen von Wirkingsfunktionen}, BMwF-GMD-22 Gesell- schaft f\"ur Mathematik und Datenverarbeitung, Bonn (1969).

\bibitem{ruzsa-graph}
I.~Z.~Ruzsa, \emph{An application of graph theory to additive number theory}, Scientia, Ser. A. \textbf{3} (1989), 97--109.

\bibitem{ruzsa-freiman}
\bysame, \emph{Generalized arithmetical progressions and sumsets}, 
Acta Math. Hungar. \textbf{65} (1994), no. 4, 379--388.

\bibitem{ruzsa-frei} \bysame, \emph{An analog of Freiman's theorem in groups}, Structure theory of set addition, Ast\'erisque \textbf{258} (1999), 323--326.

\bibitem{samorod} A.~Samorodnitsky, \emph{Low-degree tests at large distances,} available at \texttt{arXiv:math/0604353}.

\bibitem{szemeredi-3aps} E.~Szemer\'edi, \emph{Integer sets containing no arithmetic progressions,}  Acta Math. Hungar.\textbf{56}  (1990),  no. 1-2, 155--158.

\bibitem{sz-reg} \bysame, \emph{Regular partitions of graphs,} Probl\'emes combinatoires et th\'eorie des graphes, Colloq. Internat. CNRS Univ. Orsay, Orsay 1976, 399--401.

\bibitem{tao:ergodic} T.~C.~Tao, \emph{A quantitative ergodic theory proof of Szemer\'edi's theorem,} to appear in Electronic J. Combinatorics.

\bibitem{tao:icm2006} \bysame, \emph{The dichotomy between structure and randomness, arithmetic progressions, and the primes,} to appear in Proceedings of the International Congress of Mathematicians, Madrid 2006, Vol. 1.

\bibitem{tao-vu-book} T.~C.~Tao and V.~H.~Vu, \emph{Additive combinatorics,} CUP 2006.

\bibitem{ziegler}
T.~Ziegler, \emph{Universal characteristic factors and Furstenberg averages}, to appear in J. Amer. Math. Soc.


\end{thebibliography}
\end{document}